\documentclass[11pt]{amsart}
\usepackage{amsmath}
\usepackage{amssymb}
\usepackage[dvips]{graphics}
\newtheorem{Proposition}{Proposition}
  \newtheorem{Remark}[Proposition]{Remark}
  \newtheorem{Corollary}[Proposition]{Corollary}
  \newtheorem{Lemma}[Proposition]{Lemma}
  
  \newtheorem{Theorem}{Theorem}
\newtheorem{Condition}[Proposition]{Condition}
\newtheorem{Definition}[Proposition]{Definition}
\newtheorem{Assumptions}[Proposition]{Assumptions}
\newtheorem{Note}[Proposition]{Note}
\def\z{\noindent}
\def\cal{\mathcal}

\def\Box{{\hfill\hbox{\enspace${\sqre}$}} \smallskip}
\def\sqr#1#2{{\vcenter{\vbox{\hrule height .#2pt
                             \hbox{\vrule width .#2pt height#1pt \kern#1pt
                                   \vrule width .#2pt}
                             \hrule height .#2pt}}}}
\def\sqre{\mathchoice\sqr54\sqr54\sqr{4.1}3\sqr{3.5}3}
\def\mb{\mathbf}
\def\NN{\mathbb{N}}
\def\CC{\mathbb{C}}
\def\RR{\mathbb{R}}

\begin{document}
\author{ O.  Costin} \address[O. Costin]{Mathematics Department, Ohio
  State University, 231 W 18th Ave, Columbus 43210  USA}
\email{costin@math.ohio-state.edu}
\author{S. Tanveer} \title{Nonlinear evolution PDEs
  in $\RR^+\times\CC^d$: existence and uniqueness of solutions, asymptotic and
  Borel summability properties}
\address[S. Tanveer] {Mathematics Department, Ohio
  State University, 231 W 18th Ave, Columbus 43210}
\email{tanveer@math.ohio-state.edu}
\gdef\shorttitle{Nonlinear evolution PDEs in $\RR^+\times\CC^d$}
\gdef\shortauthors{O. Costin and S. Tanveer}
\maketitle \date{}

\begin{abstract}
  We consider a system of $n$-th order nonlinear quasilinear
  partial differential equations of the form
  $$
  {\bf u}_t + \mathcal{P}(\partial_{\bf x}^{\bf j}){\bf u}+{\bf g}
  \left ( {\bf x}, t, \{\partial_{\bf x}^{{\bf j}} {\bf u}\} \right )
  =0;\ {\bf {u}}({\bf x}, 0) ={\bf {u}}_I({\bf x})$$
  with
  $\mathbf{u}\in\CC^{r}$, for $ t\in (0,T)$ and large $|{\bf x}|$ in a
  poly-sector $S$ in $\mathbb{C}^d$ ($\partial_{\bf x}^{\bf j} \equiv
  \partial_{x_1}^{j_1} \partial_{x_2}^{j_2} ...\partial_{x_d}^{j_d}$
  and $j_1+...+j_d\le n$). The principal part of the constant
  coefficient $n$-th order differential operator $\mathcal{P}$ is
  subject to a cone condition. The nonlinearity ${\bf g}$ and the
  functions $\mb u_I$ and $\mb u$ satisfy analyticity and decay
  assumptions in $S$.
  
\smallskip

   The paper shows existence and uniqueness of the solution of this
  problem and finds its asymptotic behavior for large $|\bf x|$.
 
\smallskip

Under further regularity conditions on $\mb g$ and $\mb u_I$ which
ensure the existence of a formal asymptotic series solution for large
$|\mb x|$ to the problem, we prove its Borel summability to the actual
solution $\mb u$.
  
\smallskip
 
    The structure of the nonlinearity and the complex plane setting
  preclude standard methods. We use a new approach, based on
  Borel-Laplace regularization and \'Ecalle acceleration techniques to
  control the equation.
   
\smallskip

These  results are instrumental in constructive analysis of singularity
formation in nonlinear PDEs with prescribed initial data, an application
referred to in the paper.
   
\smallskip

    In special cases motivated by applications we show how the method
  can be adapted to obtain short-time existence, uniqueness and
  asymptotic behavior for small $t$, of sectorially analytic solutions,
  without size restriction on the space variable.
\end{abstract}
\vfill\eject
\tableofcontents

\section{Introduction} 
\subsection{General considerations}\label{Genc} There are relatively few general results on
existence, uniqueness and regularity of solutions of partial
differential equations in the complex domain when the conditions of
the classical Cauchy-Kowalewski (C-K) theorem are not met.  The C-K
theorem holds for first-order analytic systems (or those equivalent to
them) with analytic non-characteristic data, and for these it
guarantees local existence and uniqueness of analytic solutions.  As
is well known, its proof requires convergence of local power series
expansions. Evolution equations with higher spatial derivatives do not
satisfy the C-K assumptions and even when formal power series
solutions exist their radius of convergence is zero.  One of the goals
of this paper is to provide a theory for existence, uniqueness and
regularity of solutions in such cases, in a relatively general
setting. The theory also applies to classes of equations of higher
order in time and sufficiently high order in space after reduction (by
well known transformations, see {\em e.g.} \cite{Treves}) to evolution
systems.

The present paper generalizes \cite{CPAM} to $d$ dimensions and
arbitrary order in the spatial variable, to $r$ dimensional dependent
variable, proves additional results about short term existence and
shows Borel summability of formal solutions.   {\em A fortiori} we
obtain results on the asymptotic character of these solutions. (In  Appendix
\S\ref{illustr}, we briefly discuss
the definition and properties of Borel summation.)

Under assumptions to allow for formal expansions for large $\bf x$, we
show that series solutions are Borel summable to actual solutions of
the PDE. For this purpose we make use of \'Ecalle acceleration
techniques.  In special cases we obtain existence and uniqueness
results for $t$ in a compact set and large enough $\mathbf x$, and
separately for small $t$ and fewer restrictions on $\bf x$.

\bigskip

Properties of solutions of PDEs in the complex plane, apart from their
intrinsic interest, are relevant for properties in the real domain, as
initial singularities in $\CC$ may give rise to blow-up at later times
in the physical domain. Representation of solutions as Borel sums is
instrumental in extending techniques originally developed for ODEs
\cite{Inventiones} to find the location and type of singularities of
solutions to nonlinear PDEs \cite{CPAM3}.

It is certainly difficult to give justice to the existing theory of
nonlinear PDEs, and we  mention a number of results in the
literature relevant to the current paper. For certain classes of
PDEs in the complex domain Sammartino and Caflisch \cite{[CSI]},
\cite{[CSII]} proved the existence of nonlinear Prandtl boundary layer
solutions for analytic initial data in a half-plane.  This work
involves inversion of the heat operator $\partial_t - \partial_{YY}$
and uses the abstract Cauchy-Kowalewski theorem for the resulting
integral equation. While their method is likely to be generalizable to
certain higher-order partial differential equations, it appears
unsuitable for problems where the highest derivative terms appear in a
nonlinear manner. Such terms cannot be controlled by inversion of a
linear operator and estimates of the kernel, as used in (\cite{[CSI]},
\cite{[CSII]}).

The complex plane setting, as well as the type of nonlinearity allowed
in our paper, do not allow for an adaptation of classical, Sobolev space
based, techniques. This can be also seen in simple examples which show
that existence fails outside the domain of validity of the results we
obtain.

Certainly, many evolution equations are amenable to our setting; to
illustrate canonical form transformations and the general results we
chose a third order equation with quartic nonlinearity arising in
fluid dynamics.  Detailed singularity study \cite{CPAM3} of solutions
of this equation  relies on the present analysis.

\bigskip

Our approach extends Borel transform regularization to a general class
of nonlinear partial differential equations.  A vast literature has
emerged recently in Borel summability theory, starting with the
fundamental contributions of \'Ecalle (see e.g.  \cite{EcalleNato})
whose consequences are far from being fully explored and it is
impossible to give a quick account of the breadth of this field.  See
for example \cite{Inventiones} for more references. Yet, in the
context of relatively general PDEs, very little is known.  For small
variables, Borel summability has been recently shown for the heat
equation \cite{Lutz,Balser3}, and generalized to linear PDEs with
constant coefficients by Balser \cite{Balser2}. One large space
variable was considered by us in \cite{CPAM}, in special classes of
higher order nonlinear PDEs.  The methods in the present paper are
different and apply, for large $|\bf x|$, to a wide class of
equations.

  \subsection{Notation}\label{Nott} 
We use the following conventions. For vectors in $\CC^d$ or
multiindices we write
$$ | {\bf u} | = \sum_{j=1}^d|u_i|$$ and for multiindices we define
$$\mathbf{k}\succ \mathbf{m}\text{ if $k_i>m_i$ for all $i$}$$
If $a$
is a scalar we write $\mathbf x^a=(x_1^a,x_2^a,...,x_d^a)$. 

\z With
$\mathbf{p}$, ${\bf x}$ and ${\bf j}$ vectors of same dimension $d$,
we define
$$
{\bf p}^{\bf j} = \prod_{i=1}^dp_i^{j_i} $$
and
$$
{\partial}_{\bf x}^{\bf j} = \partial_{x_1}^{j_1}
\partial_{x_2}^{j_2}..\partial_{x_d}^{j_d} $$
We write ${\bf 1} = (1,
1, .., 1)$ and more generally, if $\alpha$ is a scalar, we write
$\boldsymbol{\alpha}= \alpha\mb 1$; thus ${\bf x^1}=\prod_{i=1}^d
x_i$. For $d$-dimensional vectors ${\bf a}$ and ${\bf b}$ we write
$$
{ \int_{\bf a}^{\bf b}} \,\,\cdot\,{\mathrm d\mb p} =
\int_{a_1}^{b_1} \int_{a_2}^{b_2} ...\int_{a_d}^{b_d} \,\,\cdot\, dp_1
dp_2 \cdots dp_d $$
The {\em directional Laplace transform} along the
ray $\arg \,p_i=\varphi_i,i=1...d$ of $F$ is given by
\begin{equation}\label{10}\left \{\mathcal{L}_{\boldsymbol \varphi} { F} \right \} ({\bf x})
\equiv \int_{\bf 0}^{{\boldsymbol \infty} e^{i {\boldsymbol \varphi}}} 
{ F} ({\bf p}) e^{-{\bf p} \cdot {\bf x}} {\mathrm d\mb p}
\end{equation}
where $\mb x e^{i\boldsymbol\theta}$ will denote the vector with
components $x_i e^{i\theta_i}$. {\em Convolution} is defined as
\begin{equation}
  \label{eq:defconv}
  (f*g)({\bf p}):=\int_{\mathbf 0}^{\bf p} 
f({\bf s})g({\bf p}-{\bf s}) {\mathrm d\bf s} 
\end{equation}
and $\sideset{^*}{}\prod$ denotes convolution product (see also
\cite{DMJ}).  Whenever used as {\em sum or product indices}, $l$ takes
all integer values between $1$ and $m$, $i$ is between $1$ and $d$, As
a sum or product multiindex, $|\mathbf j|$ indicates all $\bf j$ with
positive integer components subject to the constraint $1\le |\mathbf
j|\le n$.

\section{Problem statement and main results}\label{MainR}
\subsection{Setting and assumptions}\label{MainR1} Consider the initial value problem for a  quasilinear system
\begin{equation}\label{ufirst}
  {\bf u}_t + \mathcal{P}(\partial_{\bf x}^{\bf j}){\bf u}+{\bf g}
  \left ( {\bf x}, t, \{\partial_{\bf x}^{{\bf j}} {\bf u}\}_{|{\bf
  j}|\le n}\right ) =0; \ {\bf {u}}({\bf x}, 0) ={\bf {u}}_I({\bf
  x})\end{equation} In (\ref{ufirst}), $\mathcal{P} (\partial_{\bf x}
  ) {\bf u}$ collects the constant coefficient linear terms of the partial differential
  equation.

Emphasizing  quasilinearity, we  rewrite the equation as
\begin{multline}\label{1}
\partial_t {\bf u} + \mathcal{P} (\partial_{\bf x}) {\bf u} + 
\sum_{|{\bf J}| = n} {\bf g}_{2,{\bf J}} \left ( {\bf x}, t,
\{\partial_{\bf x}^{{\bf j}} 
{\bf u}\}_{|{\bf j}| < n} \right ) \partial_{\bf x}^{\bf J} {\bf u}\\ =
{\bf g}_1 \left ( {\bf x}, t, \{\partial_{\bf x}^{j} 
{\bf u} \}_{|{\bf j}| < n} \right);
 \ {\bf {u}}({\bf x}, 0) ={\bf {u}}_I({\bf x})
\end{multline} 

  The restrictions on ${\bf g}_1 $, ${\bf g}_2$, and ${\bf {u}_I}$
  are simpler in a normalized form, more suitable for our analysis. By
  applying $\partial_{\bf x}^{\bf j}$ to (\ref{1}) for all ${\bf j}$
  with $1\le |{\bf j}| \le n-1$, we get an extended system of
  equations for $\mb f\in\CC^m$, consisting in ${\bf {u}} $ and its
  spatial derivatives of order less than $n$, of the type (see
  Appendix for further details):
\begin{equation}\label{1.a}
\partial_t{\bf f} +\mathcal{P} (\partial_{\bf x} ) {\bf f} = 
{\sum_{{\bf q} \succeq 0 }}^\prime {\bf b_{\bf q}}
({\bf x},t,{\bf f}) \prod_{l,|{\bf j}|} 
\left (\partial_{\bf x}^{\bf j} f_l
\right )^{q_{l,{\bf j}}} + 
{\bf r} ({\bf x}, t)\ \ \mbox{with}\ \ {\bf f}({\bf x}, 0) 
={\bf f}_I ({\bf x})
\end{equation} 
where $\sum^\prime$ means the sum over the multiindices $\bf q$ with
\begin{equation}
  \label{eq:cond1}
\sum_{l=1}^m\sum_{1\le|{\bf j}|\le n} |{\bf j}| q_{l,{\bf j}} \le n  
\end{equation}
The matrix $\mathcal{P}$ is assumed to be diagonalizable, and modulo
simple changes of variables we assume it is presented in diagonal
form, $\mathcal{P}={\rm diag}\,\mathcal{P}_j,j=1,...,m$. In
(\ref{1.a}), ${\bf q} = \left (q_{l,{\bf j}} \right)$, $1\le |{\bf
  j}|\le n,1\le l\le m$ is a vector of integers and $\mathcal{P}_j$ is
an $n$-th order polynomial. We let $\mathcal{P}_{n;j}$ be the
principal part of $\mathcal{P}_j$, {\em i.e.} the part that contains
all monomials of (total) degree $n$.  The inequality (\ref{eq:cond1})
implies in particular that none of the $q_{l,{\bf j}}$ can exceed $n$
and that the summation in (\ref{1.a}) involves {\em only finitely many
  terms}.  The fact that (\ref{eq:cond1}) can always be ensured leads
to important simplifications in the proofs. Let $\rho>\rho_0>0$,
$\phi<\frac{\pi}{2n}$, $\epsilon>0$ and
\begin{equation} \label{eq:defD1} {\cal D}_{\phi, \rho;\bf x}
= \left \{\mathbf x : |\arg x_i| <\frac{\pi}{2}+\phi;\, 
|x_i| > \rho; \ i\le d\right \} 
\end{equation}
\begin{equation}\label{eq:defD}
{\cal D}_{\phi,\rho} =  {\cal D}_{\phi, \rho;\bf x}\times [0,T]
\end{equation}
\begin{Assumptions}\label{cds}

\begin{enumerate}
\item{} There is a $\phi \in \left (0, \frac{\pi}{2n} \right )$ 
such that for all ${\bf p}\ne 0$ with $\max_i|\arg p_i|\le \phi$ we have
 \begin{eqnarray}
    \label{condnew}
    \Re\, \mathcal{P}_{n;j} (-{\bf p} )>0
  \end{eqnarray}
  \item{} The functions ${\bf b}_{\bf q}(\cdot,t,\cdot) $ are analytic in $
  {\cal D}_{\frac{\pi}{2n},\rho_0}\times \{\mathbf f:|\mathbf f| <\epsilon\}$. 
We write
\begin{equation}\label{3} {\bf b}_{\bf q} ({\bf x}, t; {\bf f})
=\sum_{{\bf k} \succeq 0} {\bf b}_{{\bf q}, {\bf k}} ({\bf x},t) 
{\bf f}^{\bf k}
\end{equation}

\item{} For some constants $\alpha_r \ge 1$ independent of $T$ (see
  also \S\ref{Asympts}), $A_r(T)>0$, $\alpha_{\bf q}>0
  $\footnote{A restriction of the form $|{\bf x}|^{\tilde
      \alpha} |{\bf r} ({\bf x}, t)| < A_r (T) (*)$ may appear more
    natural. However, since every component of ${\bf x}$ is bounded
    below in $\mathcal{D}_{\phi, \rho_0, {\bf x}}$, it is clear that (*) implies
    (\ref{2}) with $\alpha_r = {\tilde \alpha}/d$.  The same comment
    applies for condition (\ref{4}). This form is more convenient in
    the present analysis. See also Note 4 following Theorem \ref{T1}.}
\begin{align}\label{2}
\sup_{\mathbf x\in {\cal D}_{\frac{\pi}{2n},\rho_0;\mathbf x} }
\left | {\bf x}^{\boldsymbol \alpha_r} {\bf r} ({\bf x}, t) \right | 
= A_r (T)<\infty \\
\sup_{\mathbf x\in {\cal D}_{\frac{\pi}{2n},\rho_0;\mathbf x} }
\left | {\bf x}^{{\boldsymbol \alpha_r}}{\bf f}_I ({\bf x},t) \right | = 
A_f (T)<\infty\label{5.1}
\\
 \sup_{\mathbf k,\mathbf q;\mathbf x\in {\cal D}_{\frac{\pi}{2n},\rho_0;
\mathbf x}}\left | {\bf x}^{{\boldsymbol \alpha}_{\bf q}} 
{\bf b}_{{\bf
q},{\bf k}} \right | = A_b (T)<\infty\label{4}
\end{align}
\item{} The analysis is interesting for $n>1$, which we assume is the
  case.

\end{enumerate}

\end{Assumptions}
\subsection{Existence and uniqueness for large $|\bf x|$}\label{EUn}
\begin{Theorem}
 \label{T1} Under the Assumptions~\ref{cds}, there is a unique solution 
 $\mathbf f$ of (\ref{1.a}) satisfying the following properties in $
 {\cal D}_{\phi,\rho_0;\mathbf x}$: (a) $\mathbf f$ analytic and (b)
 $|\mathbf x^{\bf 1}||\mathbf f|$ bounded .  Furthermore, this
 solution satisfies ${\bf f}=O({\bf x}^{-{\boldsymbol \alpha_r}})$ as
 ${\bf x} \rightarrow \infty$ in ${\cal D}_{\phi, {\tilde \rho}; {\bf
     x}} $, for large $\tilde{\rho}$.
\end{Theorem}

\noindent {\bf Notes}. {\bf 1. The proof of Theorem \ref{T1} is given in
\S \ref{Sol}} 

1. As shown in \cite{CPAM}, \cite{CPAM3} for special examples, ${\bf
  f}$, in a larger sector is expected to have singularities with an
accumulation point at infinity.

2. In section \ref{smalltime}, we also show that in some special
cases, there is a duality between {\em small} $t$ and {\em large}
${\bf x}$. 

3. Relatively simple examples in which the assumptions apply after
suitable transformations are the modified Harry-Dym equation
$H_t+H_x= H^3H_{xxx}-H^3/2$, Kuramoto-Sivashinsky $u_t + u u_x +
u_{xx} + u_{xxxx} = 0$ and thin-film equation $ h_t + \nabla \cdot
\left (h^3 \nabla \Delta h \right ) =0$ (the latter with initial
conditions such as $h({\bf x}, 0) = 1+ (1+a x_1^2 + b x_2^2)^{-1}$ in
$d=2$). The former equation is discussed in detail in \cite{CPAM} and
the normalizing process, adapted to short time analysis, is described
in \S\ref{smalltime}.

4. The condition $\alpha_r \ge 1$ is not particularly restrictive in
  problems with algebraically decaying coefficients. For these, as discussed
  in \cite{CPAM}, one can redefine ${\bf f}$ by subtracting out from it the
  first few terms of its formal asymptotic expansion for large ${\bf x}$. The
  new ${\bf f}$ decays faster at $\infty$ and the condition to $\alpha_r \ge
  1$ can be ensured.

\subsection{Borel summability of power series solutions and their asymptotic character}\label{Bsum}  Determining asymptotic properties of solutions of PDEs is
substantially more difficult than the corresponding question for ODEs.
Borel-Laplace techniques however provide a well suited modality to
overcome this difficulty.  The paper shows that formal series
solutions are  Borel summable to actual solutions (a fortiori
are asymptotic to them). A few notes on Borel summability are found in
\S\ref{illustr}.

In addition to hypothesis of Theorem \ref{T1}
we need, first of all, to impose restrictions to ensure that there exist
series solutions, to which end the coefficients of the equation should
be expandable for large $\bf x$.  In many practical applications these
coefficients turn out to be finite combinations of ramified inverse
powers of $x_i$.

\begin{Condition}\label{Cond 2}
  For large $\bf|x|$ and some $\mb N\in\NN^d$, the functions ${\bf b}_{{\bf
      q}, {\bf k}} ({\bf x}, t)$ and ${\bf r} ({\bf x}, t)$ are analytic in
  $(x_1^{-1/{N_1}},...,x_d^{-1/{N_d}})$
\end{Condition}

\begin{Theorem}\label{TrB}
  If Condition~\ref{Cond 2} and the assumptions of Theorem~\ref{T1} are
  satisfied, then the unique solution $\mb f$ found there is the Borel sum of
  its own asymptotic series. More precisely, $\mb f$ can be written as
  \begin{equation}
    \label{acc3}
    \mathbf{f}(\mathbf{x},t)=
\int_{{\RR^+}^d}e^{-\mathbf{p}\cdot\mathbf{x}^{\frac{n}{n-1}} }{\mathbf{F}_1}(\mathbf{p},t)d\mathbf{p}
  \end{equation}
  where ${\mathbf{F}_1}$ is (a) analytic at zero in
  $(p_1^{\frac{1}{nN_1}},...,p_d^{\frac{1}{nN_d}})$; (b) analytic in $\bf p\ne 0$ in the
  poly-sector $|\arg p_i|<\frac{n}{n-1}\phi+\frac{\pi}{2(n-1)}$, $i\le
  d$; and (c) exponentially bounded in the latter poly-sector.
\end{Theorem}

{\bf Comment:}  For PDEs it is known that it difficult to show, by
  classical methods, the existence of actual solutions given formal ones, when
  the formal solutions diverge.  Borel summability of a formal asymptotic
  series solution shows in particular, using Watson's lemma
  \cite{benderorszag}, that there always indeed exist {\em actual} solutions
  of the PDE asymptotic to it.  Borel summability also entails uniqueness of
  the actual solution if a sufficiently large sector of asymptoticity is
  prescribed (see, e.g., \cite{Balser}). The Borel summability parameters
  proven in the present paper are optimal, as explained in the following
  remarks, and the sharp Gevrey class of the formal solutions follows too.
\begin{Remark}
  
(i)  It follows from the same proof that $\mathbf{x}^{\frac{n}{n-1}}$ can
  be replaced with $\mathbf{x}^\beta$ for any $\beta\in
  [1,\frac{n}{n-1}]$. The canonical variable in Borel summation is
  that in which the generic Gevrey class of the formal series solution
  is one (i.e., the series diverge factorially, with factorial power
  one; \cite{Balser}).  This variable, in our case, is $\mathbf
  x^{\frac{n}{n-1}}$. 
  
  (ii) At least in simple examples, the sector of summability is
  optimal. See also Note~\ref{Opti2}.

  (iii) In many problems of interest the conditions of Theorem~\ref{TrB} are
  met by the equation in more than one sector (after suitable rotation of
  coordinates).  Then the functions $\mb F_1$ obtained in (\ref{TrB}) are
  analytic continuations of each-other, as it follows from their construction.

  (iv) If we had made the change of variable ${\bf x} \rightarrow {\bf
  x}^{n/(n-1)}$ first, (yielding the normalized  Borel variable), the
  transformed PDE would have been more difficult to handle.  Borel
  transforming directly from the $\bf x$ to $\bf p$ instead requires us
  to perform, in the proof of Theorem~\ref{TrB}, an acceleration in the sense of \'Ecalle to
  establish Borel summability, but is technically simpler.

\end{Remark}
The proof of Theorem~\ref{TrB} is given in \S\ref{B}.

\z See also \S\ref{Asympts}.  \subsection{Spontaneous formation of
    singularities in nonlinear PDEs} Borel summability of formal solutions
  associated to solutions with prescribed initial data is a key ingredient in
  the detailed analysis of spontaneous singularities of solutions and in the
  study of their global properties.  Applications of the present techniques in
  these directions, partly relying on extensions to PDEs of the methods in
  \cite{Inventiones}, are discussed in the paper \cite{CPAM3}.

\section{Inverse Laplace transform and associated integral equation}\label{S3}

The inverse Laplace transform (ILT) ${\bf G}({\bf p},t)$ of a function ${\bf
  g}({\bf x},t)$ analytic in ${\bf x}$ in ${\cal D}_{\phi, \rho;\mathbf x} $
and vanishing algebraically as ${\bf x} \rightarrow \infty$ (cf.
Lemma~\ref{L1} below and Note following it) is given by:
\begin{equation}\label{6}
{\bf G}({\bf p},t) =\left [ \mathcal{L}^{-1} \{{\bf g}\}\right ] ({\bf p}, t) 
\equiv \frac{1 }{(2 \pi i)^d } \int_{{\cal C}^d_D} 
e^{{\bf p} \cdot  {\bf x} } {\bf g}({\bf x}, t) d{\bf x} 
\end{equation}
with a contour $\mathcal{C}_D$ as in Fig. 1 (modulo homotopies),
${\cal C}_D^d \subset {\cal D}_{\phi, \rho;{\bf x}}$, and ${\bf p}$
restricted to the dual (polar) domain ${\cal S}_\phi$ defined by
\begin{equation}
  \label{eq:defSp}
  {\cal S}_\phi \equiv  \left \{ {\bf p} : |p_i|>0;\ \arg p_i \in (-\phi, \phi), 
 i=1,...,d\right \} 
\end{equation}
to ensure convergence of the integral.

The following lemma connects the $\bf p$ behavior of the ILT of 
functions of the type considered in this paper to their assumed behavior
in $\bf x$.
\begin{Lemma}
  \label{L1}
If ${\bf g}({\bf x},t)$ is analytic for ${\bf x}$ in ${\cal
D}_{\phi, \rho;\bf x}$, and satisfies
\begin{equation}\label{6.1} 
|{\bf x}^{\boldsymbol \alpha}|~  | {\bf g}(\mathbf x,t) | \le A (T) 
\end{equation}
for $\alpha\ge \alpha_0>0$, then for any $\delta\in(0,\phi)$, the ILT
$\mathbf G=\mathcal{L}^{-1}\mathbf g$ exists in $ {\cal
  S}_{\phi-\delta}$ and satisfies
\begin{equation}\label{7} 
|{\bf G}({\bf p},t)| \le C \frac{A(T)}{[\Gamma(\alpha)]^d}
|{\bf p}^{{\boldsymbol \alpha} - {\boldsymbol 1}} | e^{2 |{\bf p}|\rho}
\end{equation}
for some $C=C(\delta,\alpha_0)$.
\end{Lemma}

\begin{proof}
  The proof is a higher dimensional version of that of Lemma 3.1 in
  \cite{CPAM}.  We first consider the case when $2 \ge\alpha \ge
  \alpha_0$.  Let $C_{\rho_1}$ be a contour so that the integration
  path in each ${\bf x}$ component is as shown in Fig.~1: it passes
  through point $\rho_1+|p_i|^{-1} $, and $s=\rho_1 + |p_i|^{-1} +i
  r\exp(i\phi \mbox{\,signum($r$)})$ with $r\in(-\infty,\infty)$.
  Choosing $2\rho \ge \rho_1 \ge (2/\sqrt{3})\rho$, we have $|s|>\rho$
  along the contour and therefore, with $\arg(p_i)=\theta
  \in(-\phi+\delta,\phi-\delta)$,

$$|\mathbf g({\bf s},t)|\le A (T) |{\bf s}^{-{\boldsymbol \alpha}} | \mbox{~~and~~}
|e^{{\bf s} \cdot {\bf p}}|\le
  e^{\rho_1|\mb p|+d} e^{-r |{\bf p}|\sin|\phi+\theta|} $$ 
Thus, 
$$
\left|\int_{C_{\rho_1}} e^{{\bf s} \cdot {\bf p}}
\mathbf g({\bf s},t) {\mathrm d\bf s} \right|\le
  2 A(T) e^{\rho_1|{\bf p}|+d}\prod_{i}
   \int_0^{\infty} \left|\rho_1+|p_i|^{-1} + i r
    e^{i\phi}\right|^{-\alpha} e^{- |p_i | r\sin\delta}dr $$
\begin{equation}
\label{7.1}
\le{\tilde K} A(T) e^{\rho_1 |{\bf p}|}      
        \prod_{i} \left \{ | \rho_1 +|p_i|^{-1} |^{-\alpha} 
       \int_0^{\infty} e^{-|p_i| r \sin\delta  } dr \right \} 
\le K \delta^{-d} |{\bf p}^{{\boldsymbol \alpha} - {\boldsymbol 1}} | e^{2 \rho |{\bf p}|} 
\end{equation} 

\z where ${\tilde K}$ and $K$ are constants independent of any
parameter. Thus, the Lemma follows for $2 \ge\alpha \ge \alpha_0$, if we note
that $\Gamma (\alpha)$ is bounded in this range
of $\alpha$,  the bound only depending on $\alpha_0$.

For $\alpha >2$, there exists an integer $k >0$ so that $\alpha-k
\in (1, 2]$.  Taking 

$$[(k-1)!]^d \mathbf h({\bf x},t)=
\int_{\boldsymbol \infty}^{\bf x} 
\mathbf g({\bf z},t)({\bf x}-{\bf z})^{{\bf k}-{\bf 1}} d{\bf z}$$

\z (clearly $\mathbf h$ is analytic in ${\bf x}$, in ${\cal D}_{\phi, \rho}$
and $\partial_{\bf x}^{\bf k} {\bf h} ({\bf x}, t) = \mathbf g({\bf x},t)$), we get

\begin{multline*}
  \mathbf h({\bf x}, t) = 
\frac{(-1)^{d k} \mathbf{x}^{k\bf 1}}{[(k-1)!]^d} 
\int_{\boldsymbol 1}^{\boldsymbol \infty} \mathbf g({\bf x\cdot y}, t)
({\bf y}-{\bf 1})^{(k-1)\mb 1} d{\bf y} \\
=
\frac{(-1)^{d k}\mathbf{x} ^{(k-\alpha)\bf 1}}{[(k-1)!]^d}
\int_{\boldsymbol 1}^{\boldsymbol \infty} 
\mathbf A({\bf x\cdot y} ,t) 
{\bf y}^{-{\boldsymbol \alpha}} 
({\bf y}-{\boldsymbol 1})^{(k-1)\mb 1} d{\bf y}
\end{multline*}

\noindent with $|\mathbf A({\bf x\cdot p}, t)|\le A(T)$, whence

$$ |\mathbf h({\bf x}, t)| \le \frac{A(T) [\Gamma (\alpha-k)]^d}{
|\mathbf x^1|^{\alpha-k}
[\Gamma
(\alpha)]^d }$$

\noindent From the arguments above
with $\alpha-k$ playing the role of $\alpha$, we get
$$ |{\mathcal{L}^{-1}} \{\mathbf h\} ({\bf p}, t)| 
\le C(\delta) \frac{A(T)}{[\Gamma(\alpha)]^d}
|{\bf p^1} |^{\alpha-k-1} e^{2|{\bf p}|\rho} $$
Since ${\bf G} ({\bf p},t) =(-1)^{k d} {\mb p^{\mb 1 k}}
{{\mathcal{L}^{-1}}}\{\mathbf h\} (p, t)$,
by multiplying the above equation by $|{\bf p^1}|^{k} $, the Lemma follows for
$\alpha >2$ as well. \end{proof}

\begin{Remark}
\label{rem1}
The constant $2\rho$ in the exponential bound can
be lowered to  $\rho+0$, but (\ref{7}) suffices for 
our purposes. Note also  that the statement also holds for $\rho =0$, a fact
that will be used in \S 6. 
\end{Remark}
\begin{Remark}
\label{rem2}
Corollary~\ref{C1} 
below implies that for any ${\bf p} \in {\cal S}_\phi$, the ILT exists for the
functions ${\bf r} ({\bf x}, t)$, ${\bf b}_{{\bf q},{\bf k}} ({\bf x}, t)$, as well
as for the solution ${\bf f}({\bf x}, t)$, whose existence is shown in the
sequel.
\end{Remark}

\begin{Remark}
\label{rem3}
Conversely, if ${\bf G}({\bf p}, t)$ is any integrable
function satisfying the exponential bound in (\ref{7}), it is clear that
the Laplace Transform along a ray (\ref{10})
exists and defines an analytic function of ${\bf x}$ 
in the half-plane for each component defined by $\Re [e^{i \theta_i} x_i]
>2 \rho $ for $\theta_i \in (-\phi, \phi)$.
 Due to the width of the sector it is easy to see,
by Fubini, that $\mathcal{L}\mb G=\mb g$.
\end{Remark}

\begin{Remark}
\label{rem4}
The next corollary finds bounds for ${\bf B}_{{\bf q},{\bf
    k}}=\mathcal{L}^{-1}\{{\bf b}_{{\bf q},{\bf k}}\}$ and ${\bf
  R}=\mathcal{L}^{-1}\{{\bf r} \}$ independent of $\arg ~p_i$ for
${\bf p} \in {\cal S}_\phi$, following from the properties of
  $\bf b_{q,k} $ and ${\bf r} $ in ${\cal D}_{\frac{\pi}{2n},
    \rho_0}\supset {\cal D}_{\phi, \rho}$.
\end{Remark}

\begin{Corollary}
  \label{C1}
  The ILT of the coefficients ${\bf b}_{{\bf q},{\bf k}} $ (cf.
  (\ref{3})) and of the inhomogeneous term ${\bf r}({\bf x}, t)$ satisfy
  the following upper bounds for any ${\bf p} \in {\cal S}_\phi$

\begin{equation}\label{11}
|{\bf B}_{{\bf q},{\bf k}} ({\bf p}, t)| 
\le \frac{C_1 (\phi,\alpha_{\bf q})}{[\Gamma(\alpha_{\bf q})]^d} 
A_b (T)  |{\bf p}^{{\boldsymbol \alpha}_{\bf q} - {\bf 1} } |
e^{ {{{2}}} \rho_0 |{\bf p}|}
\end{equation}

\begin{equation}\label{12}
|{\bf R} (p, t)| \le
\frac{C_2(\phi)}{[\Gamma(\alpha_r)]^d} A_r (T)  |{\bf p}^{{\boldsymbol \alpha}_r - {\bf 1}}| 
e^{{{2}} \rho_0 |{\bf p}|}
\end{equation}
\end{Corollary}
\begin{proof}
  The proof is similar to that of Corollary 3.2 in\cite{CPAM}. From
  the conditions assumed we see that $\bf b_{q,k}$ is analytic in
  ${\bf x} \in {\cal D}_{\phi_1, \rho_0;{\bf x}}$ for any $\phi_1$
  satisfying $(2n)^{-1}\pi> \phi_1 >\phi>0$.  So Lemma~\ref{L1} can be
  applied, with $\mathbf g({\bf x},t)=\bf b_{q,k} $, with $\phi_1
  =\phi + ((2n)^{-1}\pi - \phi)/2 $ replacing $\phi$, and with
  $\delta$ replaced by $\phi_1-\phi = ((2n)^{-1}\pi - \phi)/2 $.  The
  same applies to $\mathbf R({\bf p}, t)$, leading to (\ref{11}) and
  (\ref{12}).  In the latter case, since $\alpha_r \ge~1$, $\alpha_0$
  in Lemma~\ref{L1} can be chosen to be 1.  Thus, one can choose $C_2$
  to be independent of $\alpha_r$.
\end{proof}

\begin{Lemma}\label{newcond}
For some $R\in\RR^+$ and all $\mb p$ with $|\mathbf p|>R$ and $\max_{i\le d}|\arg
  p_i|\le \phi$ we have for some $C>0$
\begin{equation}
  \label{es0}
  \Re\mathcal{P}_j(-\mathbf p)>C|\mathbf p|^n
\end{equation}
\end{Lemma}
\begin{proof}
For the proof, we take $B=\{\mathbf p:|\mathbf p|=1, \max_{j\le d} |\arg p_j|\le \phi\}$ and note that
\begin{equation}
  \label{e456}
 C_0=\inf_{\stackrel{\scriptstyle \mb p\in B }{\scriptstyle 1 \le j \le m}} \Re \mathcal{P}_{n;j}(-\mb p)>0  
\end{equation}
(cf. definitions following (\ref{eq:cond1})). Indeed, if $C_0=0$, then by
continuity $\Re \mathcal{P}_{n;j}(-\bf p)$ would have a root in $B$
which is ruled out by (\ref{condnew}). The conclusion now follows,
since on a sphere of large radius $R$, $\mathcal{P}_j$ is given by
$R^n\mathcal{P}_{n;j}(-\mb p/R)+o(R^n)$.
\end{proof}

\bigskip

\z The formal inverse Laplace transform (Borel transform) of (\ref{1.a}) 
with respect to ${\bf x}$ (see also (\ref{3})) for ${\bf p} \in \mathcal{S}_\phi$ is

\begin{equation}\label{13}
\partial_t{\bf F} +\mathcal{P} (-{\bf p}) {\bf F} 
={\sum_{{\bf q} \succeq 0}}' \sum_{{\bf k} \succeq 0}
{\bf B}_{{\bf q},{\bf k}} * {\bf F}^{*{\bf k}} 
* \sideset{^*}{}\prod_{l,|{\bf j}|}
\left ( (-{\bf p})^{\bf j} F_l \right )^{*q_{l,{\bf j}}} 
+{\bf R}({\bf p},t)
\end{equation}
where $\mathbf{F}=\mathcal{L}^{-1}\mathbf{f}$. After inverting the
differential operator on the left side of (\ref{13}) with respect to
$t$, we obtain the integral equation
\begin{multline}\label{14}
{\bf F} ({\bf p},t)={\mathcal N} ({\bf F})\equiv {\bf F}_0 ({\bf p}, t) \\
+ \int_0^t e^{-\mathcal{P} (-{\bf p}) (t-\tau)} {\sum_{{\bf q} \succeq 0}}' 
\sum_{{\bf k} \succeq 0}
{\bf B}_{{\bf q},{\bf k}}({\bf p},\tau) * {\bf F}^{*{\bf k}} ({\bf p},\tau)
* \sideset{^*}{}\prod_{l,|\mathbf j|}
\left ( (-{\bf p})^{\bf j} F_l \right({\bf p},\tau) )^{*q_{l,{\bf j}}} 
d\tau
\end{multline}
where
\begin{equation}\label{15} 
{\bf F}_0 (p,t) =e^{-\mathcal{P} (-{\bf p}) t} {\bf F}_{I} ({\bf p})
+\int_{0}^t e^{- \mathcal{P} (-{\bf p}) (t-\tau)} {\bf R} (p, \tau) d\tau \ \ 
\mbox{\z and\  ${\bf F}_I=\mathcal{L}^{-1}\{ {\bf f}_I\}$}
\end{equation}

Our strategy is to reduce the problem of existence and uniqueness
of a solution of (\ref{1.a}) to the problem of existence and uniqueness of
a solution of (\ref{14}), under appropriate conditions.

\section{Solution to the associated integral equation}
\label{Sol}

To establish the existence and uniqueness in (\ref{14}) we first introduce suitable 
function spaces.

\begin{Definition}
  \label{D1} Denoting by $\overline{\mathcal{S}_\phi}$ the closure of
  $\mathcal{S}_\phi$ defined in (\ref{eq:defSp}), $\partial
  \mathcal{S}_\phi =\overline{\mathcal{S}_\phi}\setminus
  \mathcal{S}_\phi$ and $\mathcal{K}=\overline{\mathcal{S}_\phi}\times
  [0,T]$, we define for $\nu>0$ (later to be taken appropriately large)
  the norm $ \|\, \cdot\,\|_\nu$ as
\begin{equation}\label{16} 
\| {\bf G} \|_\nu 
=M_0^d \,\sup_{({\bf p},t)\in\mathcal{K}} \left (\prod_{i} (1+|p_i|^2 \right )
e^{-\nu \,|{\bf p}|} |{\bf G}({\bf p},t)|
\end{equation}
where the constant $M_0$ (about $3.76$) is defined as
\begin{equation}\label{17}
M_0=\sup_{s\ge 0} \left\{ \frac{2(1+s^2)\left( \ln (1+s^2) + s\,\arctan s
 \right)}{s(s^2+4)}\right\}
\end{equation}

\z {\bf Note}: For fixed $\mathbf F$, $\|\mathbf F\|_\nu$ is nonincreasing
in $\nu$.
\end{Definition}

\bigskip

\begin{Definition}
  \label{D2} Consider the following Banach space.
  \begin{multline}
    {\cal A}_\phi 
\,=\,\left \{ {\bf F}: {\bf F}(\cdot,t)\mbox{ \rm analytic in }{\cal S}_\phi\right.\\\left.
    \mbox{ \rm and  continuous in }\overline{\mathcal{S}_\phi}\ 
    \mbox{ \rm for $t\in[0,T]$
    s.t.}\| {\bf F} \|_\nu <\infty \right \}
  \end{multline}

\end{Definition}

\begin{Remark}
\label{rem5}
If ${\bf G} \in {\cal A}_\phi$, then $ {\bf g}({\bf x}, t) =:{\cal
  L}_{\boldsymbol \theta} \{{\bf G} \}$ exists for suitable ${\boldsymbol
  \theta}$ if $\rho \cos (\theta_i + \arg x_i) > \nu$. Furthermore, $\mathbf
g({\bf x},t)$ is analytic in ${\bf x}$, and $|{\bf x}^{\bf 1}{\bf
  g}({\bf x}, t) |$ is bounded in ${\cal D}_{\phi, \rho;{\bf x}}$.
\end{Remark}

\begin{Lemma}
  \label{L3}
For $\nu >4 \rho_0+\alpha_r $, 
 ${\bf F}_I$ in (\ref{15})
satisfies
$$
\| {\bf F}_I \|_\nu \le C(\phi) A_{f_I} (\nu/2)^{-d \alpha_r+d} $$
while ${\bf R}$ satisfies the inequality
$$ \| {\bf R} \|_\nu  \le C(\phi) A_r (T) (\nu/2)^{-d\alpha_r+d} $$         
and therefore
\begin{equation}
\label{17.1}
\| {\bf F}_0  \|_\nu \le C(\phi) A_0 (T) (\nu/2)^{-d \alpha_r + d}     
\end{equation}
\end{Lemma}

\begin{proof} This proof is similar to that of Lemma 4.4 in \cite{CPAM}.
  We use (\ref{12}), note that $\alpha_r \ge 1$ and also that for
  $\nu >4 \rho_0+\alpha_r$ we have

\begin{multline}
\label{eq:eqnub}
\sup_{|p_1| >0} \frac{|p_1|^{\alpha_r \pm1}}{\Gamma(\alpha_r)} e^{-(\nu
-{{{2}}} \rho_0 ) |p_1|} \le \frac{(\alpha_r \pm 1)^{\alpha_r
\pm1}}{\Gamma(\alpha_r)} e^{-\alpha_r\mp 1} \left ( \nu-{{{2}}} \rho_0 \right
)^{-\alpha_r \mp 1} \\\le  K \alpha_r^{1/2\pm 1} (\nu/2)^{-\alpha_r \mp 1} 
\end{multline}
 
\z where $K$ is independent of $\nu$ and $\alpha_r$.
The latter inequality follows from  Stirling's  formula for
$\Gamma (\alpha_r)$ for large $\alpha_r$.

  Using the definition of
the $\nu-$norm and the two equations above, the inequality for $\|
\mathbf R
\|_\nu$ follows.  Since $\mathbf f_I (x)$ is required to satisfy the same bounds
as $\mathbf r(x, t)$, a similar inequality 
holds for $ \|\mathbf  F_I \|_\nu $.  Now, from the relation (\ref{15}) and
the fact that $\Re \mathcal{P}_j (-{\bf p}) $ is, by Lemma~\ref{newcond},  bounded below  for 
${\bf p} \in \mathcal{S}_\phi$,
we get the following inequality, implying (\ref{17.1})
$$ |\mathbf F_0 (\mb p, t) | \le |\mathbf F_I (\mb p) | + T {\hat A}_0 (T) 
\sup_{0 \le t \le T} |R (\mb p, t) | $$
\end{proof}
It is convenient to introduce a space of sectorially analytic
functions possibly unbounded at the origin but integrable.
\begin{Definition} Let
  \label{D3}
$$ {\cal H} 
:= \left \{ {\bf H} : {\bf H}({\bf p},t) 
~{\rm analytic ~in} ~{\cal S}_\phi,    
| {\bf H} ({\bf p},t) | \le
C \left | {\bf p}^{{\boldsymbol \alpha}-{\bf 1}} \right |   e^{\rho |{\bf p}|} \right \}  
$$
($C$, $\alpha$ and $\rho$ may depend on ${\bf H}$).
\end{Definition}

\begin{Lemma}
  \label{L4} If ${\bf H} \in {\cal H}$ and 
  ${\bf F} \in {\cal A}_\phi$, then for $\nu >\rho + 4$, for any $j$,
  ${\bf H}*F_j\in{\cal A}_\phi$, and\footnote{In the following equation, $\|\,\cdot\,\|_\nu$ is
    extended naturally to functions which are only continuous in
    $\mathcal{K}$.}:
\begin{equation}\label{18} 
\|{\bf H} * F_j \|_\nu \le \big \| |{\bf H}| * |F_j| \big \|_\nu 
\le C[\Gamma (\alpha)]^d \,\,2^{d \alpha} (\nu - \rho)^{-d \alpha} \| {\bf F} \|_\nu
\end{equation}
where $C$ is independent of $\alpha$.
\end{Lemma}

\begin{proof}
  The proof is a vector adaptation of that of Lemma 4.6 in
  \cite{CPAM}. From the elementary properties of convolution, it is
  clear that $\mathbf H*F_j$ is analytic in ${\cal S}_\phi$ and 
  continuous in $\overline{\mathcal{S}_\phi}$.  Let $\theta_i = \arg
  p_i$. We have
\begin{equation*}
| \mathbf H* F_j ({\bf p}) |  \le || \mathbf H|*|F_j | ({\bf p}) | 
\le 
\int_{\prod_{i}[0,|p_i|]} 
|\mathbf H(\mathbf s e^{i \mathbf \theta})| 
|F_j(\mathbf p -\mathbf se^{i \mathbf \theta} )|d {\bf s} 
\end{equation*}
 Now
\begin{equation}
  \label{eq:eH2}
 |\mathbf  H(\mathbf s e^{i \mathbf \theta}) | 
\le C \left | {\bf s}^{{\boldsymbol \alpha}-{\bf 1}} 
\right | e^{|{\bf s}| \rho}  
\end{equation}
\z and 
\begin{multline}
\label{18.1}
\int_{\prod_{i} [0,|p_i|]}  
{\bf s}^{{\boldsymbol \alpha} -{\bf 1}} e^{|{\bf s}| \rho} |F_j({\bf p}-{\bf s} e^{i {\boldsymbol \theta}})| 
d{\bf s}    
\\\le 
\| F_j \|_\nu e^{\nu |{\bf p}|} |{\bf p}^{{\boldsymbol \alpha}} |
\prod_i \left [ \int_0^1  \frac{ s_i^{\alpha -1} 
e^{-(\nu -\rho) |p_i| s_i } }{ M_0 (1 + |p_i|^2 (1 - s_i)^2 ) } ds_i \right ]
\end{multline}    
Since $\nu-\rho
\ge 4$, we can readily use (\ref{minilemma}) in the Appendix
with $\mu = |p_i| $, $\nu$ replaced by $\nu - \rho$, $\sigma = 1$ and $m=1$ 
to conclude 
\begin{equation}\label{19.0}
|p_i|^\alpha \int_0^1  \frac{ s_i^{\alpha -1} 
e^{-(\nu -\rho) |p_i| s_i } }{ M_0 (1 + |p_i|^2 (1 - s_i)^2 ) } ds_i 
\le \frac{K \Gamma (\alpha)\,\,2^\alpha (\nu - \rho)^{-\alpha} }{M_0 (1 + |p_i|^2 )} 
\end{equation}
Therefore, from (\ref{18.1}), we obtain
\begin{equation}
\label{19.2}
\int_{\prod_i [0,|p_i|]} 
{\bf s}^{{\boldsymbol \alpha} -{\bf 1}} e^{|{\bf s}| \rho} |F_j({\bf p}-{\bf s} e^{i {\boldsymbol \theta}})| 
{\mathrm d\bf s}    
\le 
K [\Gamma(\alpha)]^d \frac{\| F_j \|_\nu e^{\nu |{\bf p}|}\,\,2^{d \alpha} |\nu - \rho|^{-d \alpha} }{ 
M_0^d \prod_i (1+|p_i|^2)}
\end{equation}    
 From this relation, (\ref{18})
follows by applying the definition of $\| \cdot \|_\nu $. \end{proof}

\begin{Remark}
\label{remNew}
Lemma~\ref{L4} holds for $\rho=0$ as well, when $\nu > 4$.
\end{Remark}

\bigskip

\begin{Corollary}
  \label{C3}
For ${\bf F} \in {\cal A}_\phi $, and $\nu >4 \rho_0 +4 $ 
we have ${\bf B}_{{\bf q},{\bf k}}*F_{l} \in {\cal A}_\phi$ and
$$ \|{\bf B}_{{\bf q},{\bf k}}* F_{l} \|_\nu
\le \big \| |{\bf B}_{{\bf q},{\bf k}}| * |{\bf F}|\big  \|_\nu  
\le K C_1 (\phi, \alpha_{\bf q}) 
~(\nu/4)^{- d \alpha_{\bf q}} A_b (T) 
~\| {\bf F} \|_\nu $$   
\end{Corollary}

\begin{proof}
 The proof follows simply by using Lemma~\ref{L4}, with ${\bf H}$ replaced
by ${\bf B}_{{\bf q},{\bf k}}$ 
and using the relations in Corollary~\ref{C1}.
\end{proof}

\bigskip

\begin{Lemma}
\label{L4.7.1}
For ${\bf F} \in {\cal A}_\phi$, with $\nu >4 \rho_0 + 4$, for any ${\bf j}$, $l$,
$$ | {\bf B}_{{\bf q}, {\bf k}} * ({\bf p}^{\bf j} F_l ) | 
\le \frac{K C_1 | {\bf p}^{\bf j} | e^{\nu |{\bf p}|} A_b (T)}{M_0^d 
\prod_i (1 + |p_i|^2 )}
\| {\bf F} \|_\nu \left (\frac{\nu}{4} \right )^{- 
d \alpha_{\bf q}} $$  
\end{Lemma}

\begin{proof}  From the definition (\ref{eq:defconv}), it readily follows that
$$ | {\bf B}_{{\bf q}, {\bf k}} * ({\bf p}^{\bf j} F_l ) | \le |{\bf p}^{\bf j}| |{\bf
B}_{{\bf q}, {\bf k}} |* |F_l| $$

\z  The rest follows from Corollary
(\ref{C3}), and the definition of $\| \cdot \|_\nu$.
\end{proof}

\begin{Lemma}
\label{4.7.5}
For ${\bf F}$, ${\bf G} \in {\cal A}_\phi $ and $j \ge 0$
\begin{equation} \label{19.7.5}
| ({\bf p}^{{\bf j}} F_{l_1}) * G_{l_2} | \le |{\bf p}^{\bf j}|\,\, 
\big |\, |{\bf F}|*|{\bf G}|\, \big |
\end{equation}
\end{Lemma}

\begin{proof}
Let ${\bf p} =(p_1 e^{i \theta_1}, p_2 e^{i \theta_2}, ..,p_d e^{i \theta_d} )$. 
Then the result follows from the inequality
\begin{multline}
| {\bf p}^{\bf j} F_{l_1}*G_{l_2} | 
= \left| \int_{\bf 0}^{\bf p} {\bf {\tilde s}}^{\bf j} 
F_{l_1} ({\tilde {\bf s}})
G_{l_2} ({\bf p}-{\tilde {\bf s}}) d{\tilde {\bf s}} \right| 
\le  
|{\bf p}^{\bf j}| \int_{\prod_i[0,|p|_i]}
| {\bf F} (e^{i\boldsymbol\theta}\mathbf s)|  
|{\bf G} (\mathbf p-e^{i\boldsymbol\theta}\mathbf s)| d\mathbf{s}
\end{multline}
\end{proof}

\begin{Corollary}
\label{4.7.8}
If ${\bf F} \in {\cal A}_\phi $, then
\begin{equation} \label{19.7.51}
\left | \sideset{^*}{}\prod_{l,|\mathbf{j}|} 
\left ( {\bf p}^{\bf j} F_l \right )^{*q_{l,{\bf j}}} \right | 
\le\prod_i|p_i|^{\sum_{l,\mathbf{|j|}} j_i q_{l,{\bf j}}}
\left | \sideset{^*}{}\prod_{l,|\mathbf{j}|}
| {\bf F} |^{*q_{l,{\bf j}}} \right | 
\end{equation}

\end{Corollary}

\begin{proof}
This follows simply from repeated application of Lemma \ref{4.7.5}.
\end{proof}

\bigskip

\begin{Lemma}
\label{L4.7.9}
For ${\bf F}$, ${\bf G}  \in {\cal A}_\phi $,
$$ \Big | |{\bf F}| * |{\bf G}| \Big |
\le \frac{e^{\nu |{\bf p}|}}{M_0^d \prod_i (1 + |p_i|^2 )} \| {\bf F} \|_\nu 
\| {\bf G} \|_\nu $$
\end{Lemma}

\begin{proof}
\begin{equation}
\Big | |{\bf F}| * |{\bf G}| \Big | 
= \left| \int_{\bf 0}^{\bf p}  |{\bf F} ({\tilde {\bf s}})|
|{\bf G} ({\bf p}-{\tilde {\bf s}})|  {\bf d{\tilde s}}\right| 
\le   \int_{\prod_i[0,|p|_i]}
| {\bf F} (e^{i\boldsymbol\theta}\mathbf s)|  
|{\bf G} (\mathbf p-e^{i\boldsymbol\theta}\mathbf s)| d\mathbf{s}
\end{equation}
\noindent Using the definition of $\| \cdot \|_\nu$, the above
expression is bounded by
$$ \frac{e^{\nu |{\bf p}|}}{M_0^{2d}}
\| {\bf F} \|_\nu \| {\bf G} \|_\nu   \prod_i 
\int_0^{|p_i|} \frac{ds_i}{(1 +s_i^2) [1+(|p_i|-s_i)^2] } 
\le {\frac{|{\bf p}^{\bf j}| e^{\nu |{\bf p}|} }{M_0^d \prod_i (1 +|p_i|^2)} } 
\| {\bf F} \|_\nu \| {\bf G} \|_\nu    
$$
The last inequality follows from the definition (\ref{17}) of $M_0$
since
$$\int_0^{|p_i|} {\frac{ds_i}{(1 +s_i^2) [1+(|p_i|-s_i)^2] }} 
=2 
{\frac{ \ln (|p_i|^2 + 1) + |p_i| \tan^{-1} |p_i| }
 {|p_i| (|p_i|^2 +4)}} $$

\end{proof}

\bigskip

\begin{Corollary}
\label{L4.8.0}
For ${\bf F}$, ${\bf G}  \in {\cal A}_\phi $, then 
$$ \big\| |{\bf F}| * |{\bf G}| \big\|_\nu \le  \| {\bf F} \|_\nu 
\| {\bf G} \|_\nu $$
\end{Corollary}

\begin{proof}
  This is an application of Lemma \ref{L4.7.9} and the definition of
  $\| \cdot \|_\nu$.
\end{proof}

\bigskip
\begin{Lemma}
\label{L4.9.5}
For $\nu >4 \rho_0 + 4 $,  
\begin{multline} \label{21.5}
\left | {\bf B}_{{\bf q}, {\bf k} } * {\bf F}^{*{\bf k}} 
* \!\!\sideset{^*}{}\prod_{l,{\bf |j|}}
\left ({\bf p}^{\bf j} F_l \right )^{*q_{l,{\bf j}}} \right | 
\le \frac{e^{\nu |{\bf p}|} \prod_i |p_i|^{\sum j_i q_{l,{\bf j}} } }{
M_0^d \prod_i (1 + |p_i|^2 ) }   
\| {\bf F} \|_\nu^{|\bf q| + |{\bf k}| - 1} 
\big\| |{\bf B}_{{\bf q},{\bf k}} | * |{\bf F} | \big\|_\nu 
\end{multline}

\z if $(\mathbf q, \mathbf k)\ne (\mathbf0, \mathbf 0)$ and is zero if
$(\mathbf q, \mathbf k)= (\mathbf0, \mathbf0)$.

\end{Lemma}

\begin{proof} For ${\bf (q,k)} ={\bf (0,0)} $  we have ${\bf B}_{{\bf q},{\bf k}}=0$ 
  (see remarks after eq. (\ref{3})).  If ${\bf k} \ne {\bf 0}$,
  Corollary \ref{4.7.8} shows that the left hand side of (\ref{21.5})
  is bounded by
 $$
 \prod_i|p_i|^{\sum j_i q_{l,{\bf j}}} \left | |{\bf B}_{{\bf
       q}, {\bf k}}|*|{\bf F}| * |{\bf F}|^{*(|{\bf
       k}|-1)}*\!\!\sideset{^*}{} \prod_{l,{\bf |j|}} |{\bf F}|^{*q_{l,{\bf j}}}
 \right | $$
 Using Corollaries \ref{4.7.8} and \ref{L4.8.0} and Lemma
 \ref{L4.7.9}, the proof follows for ${\bf k} \ne 0$.  Similar steps
 work for the case ${\bf k} = {\bf 0}$ and ${\bf q} \ne {\bf 0}$,
 except that ${\bf B}_{{\bf q},{\bf k}}$ is convolved with ${\bf
   p}^{{\bf j}'} F_{l_1}$ for some $({\bf j}', l_1)$, for which the
 corresponding $q_{l_1,{\bf j}'} \ne 0$, and we now use Lemma
 \ref{4.7.5} and the definition of $\| \cdot \|_\nu$.
\end{proof}

\bigskip
\begin{Corollary}
\label{C5}
For $\nu >4 \rho_0 + 4 $, 
\begin{multline} 
\left | {\bf B}_{{\bf q}, {\bf k} } * {\bf F}^{*{\bf k}} 
* \sideset{^*}{}\prod_{l,{\bf |j|}}
\left ({\bf p}^{\bf j} F_l \right )^{*q_{l,{\bf j}}} \right | 
\\ \le \frac{K C_1 A_b (T) e^{\nu |{\bf p}|} \prod_i 
|p_i|^{\sum j_i q_{l,{\bf j}} } }{
M_0^d \prod_i (1 + |p_i|^2 ) }  
\left ( \frac{\nu}{4} \right )^{-d \alpha_{\bf q} } 
\| {\bf F} \|_\nu^{|q| + |{\bf k}|} 
\end{multline}
\end{Corollary}
The proof follows immediately from Corollary \ref{C3} and
Lemma \ref{L4.9.5} $\Box$.

\begin{Lemma}  \label{L6}
For $\nu >4 \rho_0 + 4$, we have
\begin{multline}\label {22}
  \left| \int_{0}^t e^{-\mathcal{P} (-{\bf p}) (t-\tau)} {\bf B}_{{\bf
        q},{\bf k}} * {\bf F}^{*{\bf k}} *
    \sideset{^*}{}\prod_{l,\bf|j|}
    \left ({\bf p}^{\bf j} F_l \right )^{*q_{l,{\bf j}}} d\tau \right | \\
  \le  \frac{C
    \tilde{A}_b (T) e^{\nu |{\bf p}|} }{M_0^d \prod_i (1 +
    |p_i|^2 ) } \left ( \frac{\nu}{4} \right )^{-d
    \alpha_{\bf q}} \| {\bf F} \|_\nu^{|\boldsymbol q|+ |{\bf k}| }
\end{multline}
for some $\tilde{A}_b(T)\ge A_b(T)$ (evaluated in the proof) and where the
constant $C$ is independent of $T$, but depends on $\phi$ and
$\alpha_{\bf q}$.
\end{Lemma}

\begin{proof} This is a consequence of  Lemmas
  \ref{L4.7.1} and \ref{L4.9.5} and the fact that for $0 \le |{\bf
    l'}| \le n$ we have, for $|\mathbf p| \le R$ (with $R$ as in
  Lemma \ref{newcond}),
\begin{equation}
\label{eq:28.5.0}
 J:=
    \left | \mathbf p^{\mathbf l'} \right | \int_0^t e^{-\Re \mathcal{P} (-{\bf p})
      (t-\tau)} d \tau  \le C_2(T)
\end{equation}
For $|\mb p|>R$ we have, by Lemma~\ref{newcond},
  $\mathcal{P} (-{\bf p})>C|\bf p|^n$, and $J$ is majorized by
\begin{multline}
\label{eq:28.5}
m\max_{j\le m} \frac{ |\mathbf p^{\mathbf l'}  |}{ \Re \mathcal{P}_j (-{\bf p}) }  
\left [ 1 - e^{-\Re \mathcal{P}_j (-{\bf p}) t}  \right ]  \le 
\max_{j\le m}\frac{T^{1 - |{\bf l'}| /n} |{\bf p}|^{|\bf l'| } }{ 
|\Re \mathcal{P}_j (-{\bf p}) 
|^{|{\bf l'}|/n}} 
\sup_{\gamma > 0} \frac{1 - e^{-\gamma}}{\gamma^{1 - | {\bf l'} | /n } } 
\\\le C T^{1 - |{\bf l'}| /n}
\end{multline}
where ${\bf l'} = \sum_{{\bf j}, l} {\bf j} q_{l, {\bf j}} $.
\end{proof}

\bigskip
\begin{Definition}
  \label{D4}
For ${\bf F}$ and ${\bf h}$ 
in ${\cal A}_\phi$, and ${\bf B}_{{\bf q},{\bf k}} \in {\cal H}$, as above,
define ${\bf h}_0 = {\bf 0}$ and for $k \ge 1$, 
\begin{equation}\label {23}
{\bf h}_{\bf k} 
\equiv {\bf B}_{{\bf q},{\bf k}}
*[ ({\bf F}+{\bf h})^{*{\bf k}} - {\bf F}^{*{\bf k}}].  
\end{equation} 
\end{Definition}

\bigskip
\begin{Lemma}
  \label{L7}
For $\nu >4 \rho_0 + 4$, and for ${\bf k} \ne 0$,
\begin{equation}\label {24}
\| {\bf h}_{\bf k} \|_\nu \le  
|{\bf k}| {
\Big( \| {\bf F} \|_\nu + \| {\bf h} \|_\nu \Big) }^{|{\bf k}|-1} 
\big\| |{\bf B}_{{\bf q},{\bf k}}|*|{\bf h}| \big\|_\nu 
\end{equation}
and is zero for ${\bf k} = 0$.
\end{Lemma}
  
\begin{proof}
  The cases $|{\bf k}|=0,1$ follow from the definition of ${\bf
      h}_0$ and (\ref{23}) respectively.  Assume formula (\ref{24})
  holds for all $ |{\bf k}| \le l$. Then all multiindices of length
  $l+1$ can be expressed as ${\bf k} + {\bf \hat{e}}_i$, where ${\bf
    \hat{e}}_i\in\RR^m$ is the $m$ dimensional unit vector in the
  $i$-th direction, and $|{\bf k}|=l$.
  \begin{multline*}
     \| {\bf h}_{{\bf k} +{\bf \hat{e}}_i} \|_\nu = \| {\bf B}_{{\bf q},{\bf
k}}*(F_i+h_i) *({\bf F}+{\bf h})^{*{\bf k}} - {\bf B}_{{\bf q},{\bf k}}
*F_i*{\bf F}^{*\bf k} \|_\nu \\=\| {\bf B}_{{\bf q},{\bf k}} *h_i*({\bf
F}+{\bf h})^{*{\bf k}} + F_i*{\bf h}_{\bf k} \|_\nu  
  \end{multline*}

\z Using (\ref{24}) for $|{\bf k}|=l$, we get
$$ \le \| |{\bf B}_{{\bf q},{\bf k}} |*|{\bf h} | \|_\nu 
\left ( \| {\bf F} \|_\nu + \| {\bf h} \|_\nu \right )^l 
+ l \| {\bf F} \|_\nu \Big( \| {\bf F} \|_\nu + \|{\bf h} \|_\nu \Big)^{l-1} 
 \big\| |{\bf B}_{{\bf q},{\bf k}}|*|{\bf h}|  \big\|_\nu
$$
$$
\le 
(l + 1) {
\Big( \| {\bf F} \|_\nu + \| {\bf h} \|_\nu \Big) }^{l} 
\big\| |{\bf B}_{{\bf q},{\bf k}}|*|{\bf h}| \big\|_\nu 
$$
Thus (\ref{24}) holds for $|{\bf k}|=l+1$. 
\end{proof}

\bigskip
\begin{Definition}
  \label{D41}
  For ${\bf F}\in {\cal A}_\phi$ and ${\bf h}\in {\cal A}_\phi$, and
  ${\bf B}_{{\bf q},{\bf k}}$ as above define ${\bf g}_{\bf 0} = {\bf
    0}$, and for $ |\mathbf q| \ge 1$,
\begin{equation}\label {23.5}
{\bf g}_{\bf q} 
\equiv {\bf B}_{{\bf q},{\bf k}}* \sideset{^*}{}\prod_{l,\bf |j|} \left ({\bf p}^{\bf j} [F_l + h_l]\right )^{*q_{l,{\bf j}}} - 
{\bf B}_{{\bf q},{\bf k}}* \sideset{^*}{}\prod_{l,\bf |j|}
\left ({\bf p}^{\bf j} F_l \right )^{*q_{l,{\bf j}}} 
\end{equation} 
\end{Definition}

\bigskip

\begin{Lemma}
  \label{L7.5}
For $\nu >4 \rho_0 + 4$, ${\bf g}_{\bf 0}=0$ and for
$|{\bf q}| \ge 1 $
\begin{equation}\label {24.5}
\left | {\bf g}_{\bf q} \right | \le  
\left |{\bf p}^{\sum {\bf j} q_{l, {\bf j}} }\right |
\frac{ e^{\nu |{\bf p}|} |{\bf q}|}{M_0^d \prod_i (1+|p_i|^2)}      
{\Big( \| {\bf F} \|_\nu + \| {\bf h} \|_\nu \Big) }^{|{\bf q}|-1} 
\big\| |{\bf B}_{{\bf q},{\bf k}}|*|{\bf h}| \big\|_\nu 
\end{equation}
and is zero for ${\bf q} = 0$.
\end{Lemma}

\begin{proof}
  The cases $|{\bf q}| =0,1 $ follow from the definition of
    ${\bf g}_{\bf 0}$ and (\ref{23.5}) respectively (since only terms
    linear in ${\bf F}$ are involved in (\ref{23.5})). Assuming
  (\ref{24.5}) holds if $|{\bf q}|\le l$ we show that it
    holds for ${\bf q} + {\bf \hat{e}} $,
    where ${\bf \hat{e}}$ is a unit vector, say in the $(l_1,
    j'_1,j'_2,...,j'_d)$ direction.  We have
\begin{multline}\label{ML1}
\left | {\bf g}_{{\bf q} +{\bf \hat{e}}} \right |  \le
\Bigg | 
{\bf B}_{{\bf q},{\bf k}}*\left [ {\bf p}^{{\bf j}'} (F_{l_1} + h_{l_1} ) \right ]*    
\sideset{^*}{}\prod_{l,\bf |j|}
\left [ {\bf p}^{\bf j} (F_l + h_l) \right ]^{*q_{l,{\bf j}}} 
\\-{\bf B}_{{\bf q},{\bf k}}*\left [ {\bf p}^{{\bf j}'} F_{l_1} \right ]*    
\sideset{^*}{}\prod_{l,\bf |j|}
\left [ {\bf p}^{\bf j} F_l \right ]^{*q_{l,{\bf j}}} \Bigg |\\  
\le \left | {\bf B}_{{\bf q},{\bf k}}* \left ( {\bf p}^{{\bf j}'} h_{l_1} \right
) \right | * \left |
\sideset{^*}{}\prod_{l,\bf |j|} \left [
{\bf p}^{\bf j} (F_l + h_l) \right ]^{*q_{l,{\bf j}}} \right | + | ({\bf p}^{{\bf j'}}
F_{l_1})*{\bf g}_{\bf q} | 
\end{multline}
 
\z 
Using Lemma \ref{L4.9.5} and equation (\ref{24.5}), we get the following upper bound 
implying the induction step
$$ 
\left | {\bf g}_{{\bf q} +{\bf \hat{e}}} \right |\le \frac{\left | {\bf p}^{{\bf j}' 
+ \sum {\bf j}  q_{l,{\bf j}}} \right | e^{\nu |{\bf p}|}}{M_0^d \prod_i (1+|p_i|^2)}  
\left ( \| {\bf F} \|_\nu + \| {\bf h} \|_\nu \right )^{\sum q_{l,{\bf j}}}
 \big\| |{\bf B}_{{\bf q},{\bf k}}|*|{\bf h}|  \big\|_\nu  
$$
$$
+\frac{\left |{\bf p}^{{\bf j}'+ \sum {\bf j} q_{l,{\bf j}}} \right| 
|\mathbf q| e^{\nu |{\bf p}|}}{M_0^d \prod_i (1+|p_i|^2)}  
\left ( \| {\bf F} \|_\nu + \| {\bf h} \|_\nu \right )^{|\mathbf q| - 1}
\|{\bf F}\|_\nu  
 \big\| |{\bf B}_{{\bf q},{\bf k}} |*|{\bf h}|  \big\|_\nu 
$$ 
$$
\le \frac{|{\bf p}^{\sum {\bf j} (q_{l,{\bf j}}+e_{l,{\bf j}})} (|{\bf q}+{\bf
    \hat{e}}| e^{\nu |{\bf p}|}}{\prod_i M_0^d \prod_i (1+|p_i|^2)} 
  \left ( \| {\bf F} \|_\nu +
  \| {\bf h} \|_\nu \right )^{|\bf q| } \big\| |{\bf
  B}_{{\bf q},{\bf k}}|*|{\bf h}| \big\|_\nu $$
\end{proof}

\begin{Lemma}
  \label{L8}
For ${\bf F}$ and ${\bf h}$ in ${\cal A}_\phi$, $\nu >4 \rho_0 + 4$,
\begin{multline*}
  \Bigg | {\bf B}_{{\bf q},{\bf k}}
*\Big({\bf F} + {\bf h}\Big)^{*{\bf k}}*\sideset{^*}{}\prod_{l,\bf |j|}
\left ({\bf p}^{\bf j} (F_l + h_l) \right )^{*q_{l,{\bf j}}} -
{\bf B}_{{\bf q},{\bf k}}
*{\bf F}^{*{\bf k}}*\sideset{^*}{}\prod_{l,\bf |j|}
\left ({\bf p}^{\bf j} F_l \right )^{*q_{l,{\bf j}}} \Bigg | 
\end{multline*}
\begin{equation}\label {25}
\le \frac{\left | {\bf p}^{\sum {\bf j} q_{l,{\bf j}}} \right |
(|\mathbf q|+|{\bf k}|) 
e^{\nu |{\bf p}|}}{M_0^d \prod_i (1+|p_i|^2)} 
\left ( \| {\bf F} \|_\nu
+ \| {\bf h} \|_\nu \right )^{|{\bf k}|+|\mb q| - 1} 
\| |{\bf B}_{{\bf q},{\bf k}}|*|{\bf h}| \|_\nu 
\end{equation}

\z if $({\bf q},{\bf k}) \ne (\bf 0,0)$ and is zero otherwise.
\end{Lemma}

\bigskip
\begin{proof}
It is clear from (\ref{23}) that the left side of (\ref{25}) is simply
$$ \left | {\bf h}_{\bf k} * \sideset{^*}{}\prod_{l,\bf |j|}
\left ( {\bf p}^{\bf j} (F_l + h_l ) \right )^{*q_{l,{\bf j}}} 
+{\bf F}^{*{\bf k}}*{\bf g}_{\bf q} \right | $$
However, from Corollary \ref{4.7.8}, Lemmas \ref{L4.7.9} and \ref{L7},
\begin{multline*}
 \!\!\! \!\!\! \left | {\bf h}_{\bf k} \!*\!\!\!\sideset{^*}{}\prod_{l,\bf |j|}
\left ( {\bf p}^{\bf j} (F_l + h_l ) \right )^{*q_{l,{\bf j}}} \right |
\le \frac{\left | {\bf p}^{\sum {\bf j} q_{l,{\bf j}}} \right |
|{\bf k}| e^{\nu |{\bf p}|}}{M_0^d \prod_i  (1+|p_i|^2)} 
\left ( \| {\bf F} \|_\nu
+ \| {\bf h} \|_\nu \right )^{|{\bf k}|+|\mb q| - 1} 
\big\| |{\bf B}_{{\bf q},{\bf k}}|*|{\bf h}| \big\|_\nu 
\end{multline*}
and from Corollary \ref{4.7.8}, Lemmas \ref{L4.7.9} and \ref{L7.5},

$$ \left | {\bf F}^{*{\bf k}}*{\bf g}_{\bf q} \right | \le
\frac{\left| {\bf p}^{\sum {\bf j} q_{l,{\bf j}}} \right | 
|\mb q|e^{\nu |{\bf p}|}}{M_0^d
\prod_i (1+|p_i|^2)} \left ( \| {\bf F} \|_\nu + \| {\bf h} \|_\nu \right
)^{|{\bf k}|+|\mb q| - 1} \big\| |{\bf B}_{{\bf
q},{\bf k}}|*|{\bf h} | \big\|_\nu $$

\z  Combining these two inequalities, the
proof of the lemma follows. \end{proof}

\bigskip 
\begin{Lemma} 
\label{L8.9}
For $\nu >4 \rho_0 + 4$ we have
\begin{multline}
\left\| \int_0^t e^{-\mathcal{P} (-{\bf p}) (t-\tau)} 
\left [ {\bf B}_{{\bf q},{\bf k}}
*\Big({\bf F} + {\bf h}\Big)^{*{\bf k}}*\sideset{^*}{}\prod_{l,\bf |j|}
\left ({\bf p}^{\bf j} (F_l + h_l) \right )^{*q_{l,{\bf j}}}\right.\right.\\\left.\left.-
{\bf B}_{{\bf q},{\bf k}}
*{\bf F}^{*{\bf k}}*\sideset{^*}{}\prod_{l,\bf |j|}
\left ({\bf p}^{\bf j} F_l \right )^{*q_{l,{\bf j}}} \right ]  d\tau \right\|_\nu \\   
\le \tilde{A}_b (T) 
C(\phi) (|\mb q|+|{\bf k}|) 
\left ( \| {\bf F} \|_\nu
+ \| {\bf h} \|_\nu \right )^{|{\bf k}| +|\mb q| - 1} 
\left ( \frac{\nu}{4} \right )^{ 
-d \alpha_{\bf q}} \| {\bf h} \|_\nu  
\end{multline}
\end{Lemma}

\begin{proof} This follows from  Corollary \ref{C3}  and Lemma
  \ref{L8} and the definition of $\| \cdot \|_\nu$ together with the
  bounds (\ref{eq:28.5.0}) and (\ref{eq:28.5}). \end{proof}
 \begin{Lemma}
  \label{L10}
  
  For ${\bf F} \in {\cal A}_\phi$, and $\nu >4 \rho_0 + \alpha_r+3$ large
  enough, 
(see Note after Definition
  (\ref{D1})), ${\mathcal N} ({\bf F}) $ defined in (\ref{14}) satisfies
  the following bounds
\begin{equation}\label{26.8} 
\| {\mathcal N} ({\bf F}) \|_\nu \le \|{\bf F}_0 \|_\nu
+ C(\phi) \tilde{A}_b (T)
{\sum_{{\bf q} \succeq 0}}' \sum_{{\bf k} \succeq 0}
\left ( \frac{\nu}{4} 
\right )^{-d \alpha_{\bf q}} \|\mathbf  F \|_\nu^{|\mb q|+|\mb k|}
\end{equation}

\begin{multline} \label{26.9}
\| {\mathcal N} ({\bf F} + {\bf h}) - {\mathcal N} ({\bf F}) \|_\nu
\le C(\phi) \tilde{A}_b (T)  \| {\bf h} \|_\nu  \times \\
{\sum_{{\bf q} \succeq 0}}' \sum_{{\bf k} \succeq 0} 
\left ( \frac{\nu}{4} \right )^{
- d \alpha_{\bf q} } ( |\mb q| + | {\bf k} | ) 
\left ( \| {\bf F} \|_\nu + \| {\bf h} \|_\nu \right )^{|\mb q|
+ |{\bf k}| - 1} 
\end{multline}
\end{Lemma}

\begin{proof}
  The proofs are immediate from the expression (\ref{14}) of
  ${\mathcal N} ({\bf F})$ and Lemmas \ref{L6}, \ref{L7} and
  \ref{L8.9}.  
  Note also that the sum with respect to ${\bf
    q}$ only involves finitely many terms, see (\ref{eq:cond1}). 

\end{proof}

\bigskip

\begin{Remark}
\label{rem7}
Lemma \ref{L10} is the key to showing the
existence and uniqueness of a solution in ${\cal A}_\phi$ to (\ref{14}),
since it provides the conditions for the nonlinear operator ${\mathcal
N}$ to map a ball into itself as well the necessary contractivity
condition.
\end{Remark}

\bigskip
\begin{Lemma}
  \label{L11}
If there exists some $b >1$ so that 
\begin{equation}\label{29} 
b \| \mathbf F_0 \|_\nu < 1\end{equation}
\z and
\begin{equation}\label{29.01} 
C(\phi) \tilde{A}_b (T) {\sum_{{\bf q} \succeq 0}}' \sum_{{\bf k} \succeq 0} 
\left (\frac{\nu}{4} 
\right )^{-d \alpha_{\bf q}} 
\| b {\bf F}_0 \|_{\nu}^{|{\bf k}|+|\mb q|}
<1 - \frac{1}{b} 
\end{equation}
then the nonlinear mapping ${\mathcal N}$, as defined in (\ref{14}), maps a ball of
radius $b \| {\bf F}_0 \|_\nu$ into itself. Furthermore, if 
\begin{equation} \label{29.1}
C(\phi) \tilde{A}_b (T) {\sum_{{\bf q} \succeq 0}}' \sum_{{\bf k} \succeq 0} 
(|{\bf q}|+|{\bf k}|) \left (\frac{\nu}{4} 
\right )^{-d \alpha_{\bf q}}
{(3 b)}^{|{\bf k}|+|\mb q|-1}
\| {\bf F}_0 
\|_{\nu}^{|{\bf k}|+|\mb q|-1}
<1
\end{equation}
then ${\mathcal N}$ is a contraction there.
\end{Lemma}
 
\begin{proof}
This is a simple application of Lemma~\ref{L10},
if we  note that in the ball of radius $ b \| {\bf F}_0 \|$,
$ \|{\bf F} \|_\nu^k < b^k \| {\bf F}_0 \|_\nu^k $ 
and using in (\ref{26.9}) the fact that  
$\| {\bf F} \|_\nu + \| {\bf h} \|_\nu \le 3 b \|{\bf F}_0 \|_\nu $ 
if  $\max\{\|{\bf F}\|_\nu,\|{\bf F} + {\bf h}\|_\nu\}<b\|\mb F_0\|$.
\end{proof}

\bigskip
\begin{Lemma}
  \label{L12}
Consider $T >0$ and $\phi\in (0, (2n)^{-1}\pi)$ so that
(\ref{condnew}) is satisfied. Then,
for all sufficiently
large $\nu$, there exists a unique ${\bf F} \in {\cal A}_\phi$ that
satisfies the integral equation (\ref{14}). 
\end{Lemma}

\begin{proof}
  We choose $b=2$ for definiteness.  It is clear from the bounds on
  $\|{\bf F}_0 \|_\nu$ in Lemma \ref{L3} 
that for given $T$, since $\alpha_r \ge
  1$, we have $ b \|{\bf F}_0 \|_\nu <1$ for all
  $\nu $ large. Further, it is clear by inspection that all conditions
  (\ref{29}), (\ref{29.01}) and (\ref{29.1}) are satisfied for all
  sufficiently large $\nu$. The lemma now follows from the contractive
  mapping theorem. \end{proof}

\subsection{Behavior of  $^s\!\mathbf F$ near $\mathbf p=0$}

In the following proposition, we denote by 
$^s\!\mathbf F$ the solution $\bf F$ of Lemma 
\ref{L12}. 

\begin{Proposition}
  \label{AsymptAtp=0} For some $K_1>0$ and small ${\bf p}$ we have
  $|^s\!\mathbf F|\le K_1|{\bf p^1}|^{\alpha_r-1} $ and thus 
$|^s\!\mathbf f|\le K_2
  |{\bf x^1}|^{-\alpha_r}$ for some $K_2>0$ in ${\cal D}_{\phi,\rho}$ as
  $|{\bf x}|\rightarrow\infty$.
\end{Proposition}
\begin{proof}
  The idea of the proof is to note that, once we have found $^s\!{\bf F}$, 
this function also satisfies in a neighborhood 
of the origin $\overline{\mathcal{S}}_a=\overline{\mathcal{S}}\cap \{{\bf p}:\,|p_i|\le a_i \}$ 
a linear equation of the form
\begin{equation}
\label{29.2}
^s\!{\bf F}=\mathcal{G}\, (^s\!{\bf F})+{\bf F}_0\ \ \mbox{\rm or }\ \
^s\!{\bf F}=(1-\mathcal{G})^{-1}{\bf F}_0
\end{equation}
where, of course, ${\cal G}$ depends on the previously found $^s\!{\bf
  F}$; there are many choices of ${\cal G}$ that work. Every term in
the sum in (\ref{14}) is a convolution product; in each of them we
replace all but one component of $\mathbf{F}$ by the corresponding
component of $\mathbf{^sF}$; ${\cal G} {\bf F}$ is defined as the sum
of the terms thus constructed. Estimates of the form used for
Lemma~\ref{L10} show uniform convergence of the sum for large enough
$\nu$ (or small $\bf a$). The result is a $\mathcal{G}$ as below,
where the sum over $\boldsymbol{\mu}$ contains only finitely many
terms and which has manifestly small norm if $\bf a$ is small (or
$\nu$ is large)
\begin{equation}
  \label{eq:defcalg}
   {\cal G} {\bf F} = \int_0^{t} e^{-\mathcal{P} (-{\bf p}) (t-\tau)} 
\left [ \sum_{l}
{\bf G}_l * F_l + \sum_{\boldsymbol \mu}{\bf {\hat G}}_{\boldsymbol \mu}
*\left ( (-{\bf p})^{\boldsymbol \mu} F_{l_{\boldsymbol \mu}} \right ) \right ] d\tau
\end{equation}
By (\ref{2}), (\ref{5.1}), (\ref{15}) and Lemma~\ref{L1}, we see that
$\|{\bf F}_0\|_{\infty}\le K_3 \left | {\bf a}^{\alpha_r-1} \right |$
in $\overline{\mathcal{S}}_{\bf a}$ for some $K_3>0$ independent of
${\bf a}$.  Then, from (\ref{29.2}) for small enough $|{\bf
  a}|$, we have

$$\max_{\overline{\mathcal{S}}_a}|^s\!{\bf F}(p,t)|=\|^s\!{\bf F}\|\le
(1-\|\mathcal{G}\|)^{-1}\max_{\overline{\mathcal{S}}_a}\|{\bf
  F}_0\|\le 2K_3 |{\bf a}^{\alpha_r-1} |$$
\z and thus for small
$|{\bf p}|$, we have $|{\bf F}({\bf p},t)|\le 2K_3 \left | {\bf
    p}^{\alpha_r-1} \right |$ and the proposition follows. Indeed, the
arguments also show that that the same estimates hold when any
component $p_i \rightarrow 0$, if the others are bounded.
\end{proof}

\subsection{End of proof of Theorem~\ref{T1}.} 
Lemma~\ref{L1} shows that if ${\bf f}$ is a solution of
  (\ref{1.a}) satisfying $|{\bf x}^{\bf 1} | |{\bf f}| \le A(T)$ for
${\bf x} \in \mathcal{D}_{\phi, \rho, {\bf x}}$,
then ${{\mathcal{L}^{-1}}}
  \{{\bf f}\}\in \mathcal{A}_{\phi-\delta}$ for $ 0 < \delta < \phi$
  for $\nu$ sufficiently large.  For large enough 
$\rho$, the series (\ref{3})
  converges uniformly for ${\bf x} \in \mathcal{D}_{\phi, \rho, {\bf x}}$ 
  and thus ${\bf F}={{\mathcal{L}^{-1}}} \{{\bf
    f}\}$ satisfies (\ref{14}), which by Lemma~\ref{L12} has a unique
  solution in $\mathcal{A_\phi}$ for any $\phi$ $\in$ $(0,
  ~(2n)^{-1}\pi)$ for which (\ref{condnew}) holds.  Conversely, if
  $^s\!{\bf F}\in\mathcal{A}_{\tilde \phi}$ is the solution of
  (\ref{14}) for $\nu>\nu_1$, then, for sufficiently large $\rho$,
  $^s\!{\bf f}=\mathcal{L}\,^s\!{\bf F}$ is analytic in ${\bf x}$ in
  $\mathcal{D}_{\phi,\rho}$ for $0< \phi < {\tilde \phi} <
  {(2n)^{-1}\pi} $ (cf. Remark~\ref{rem5}). Proposition
  \ref{AsymptAtp=0} shows that $^s\!\mathbf f=O({\bf x}^{-{\bf \alpha_r}
    })$ and entails uniform convergence of the series in (\ref{1.a}).
  By the properties of Laplace transforms, $^s\!{\bf f} $ solves the
  problem (\ref{1.a}).

\section{Borel summability of formal solutions to the PDE}\label{B}
We now assume Condition 1 in addition to Assumption 1.  In our
approach it was technically convenient to use oversummation, in that
the inverse Laplace transform was performed with respect to
$\mathbf{x}$.  Showing Borel summability in the appropriate variable
($\mathbf{x}^{\frac{n}{n-1}}$, as explained) requires further
arguments.
\subsection{Behavior of $\mathbf{F}$ for large $|\bf p|$ outside $\mathcal{S}_\phi$}\label{701} For the purpose of showing Borel summability of formal series solutions we need to control
$\mathbf{F}$ for large $|\bf p|$ uniformly in $\CC^d$. For this
purpose we introduce two other Banach spaces, relevant to the
properties we are aiming to show.  Firstly, let
$\mathfrak{B}(\nu,n,\mathcal{S})$ be the Banach space of functions
analytic in the sector $\mathcal{S}=\{\mathbf{p}:|p_i|>0, \arg(p_i)\in
(a_i,b_i)\}$ and continuous in its closure, where $b_i-a_i$ will be
chosen larger than $2\pi N_i$ (cf. Condition~\ref{Cond 2})
The Banach space is equipped with the norm
\begin{equation}
  \label{normnu,n}
\|\boldsymbol{\Psi}\|_{\nu n}=\sup_{\mb p\in\mathcal    S;t\in[0,T]}\left|\boldsymbol{\Psi}({\bf p},t)e^{-\nu (t+1) \sum_{j}(|p_j|+|p_j|^n)}\right|
\end{equation}

\begin{Lemma}\label{combined}
  For any intervals $(a_i, b_i), $ $ i=1,...,d$ the solution $\mathbf{F}$
  of (\ref{14}) given in Lemma ~\ref{L12} is in
  $\mathfrak{B}(\nu,n,\mathcal{S})$.
\end{Lemma}
\begin{proof}
  Because of the obvious embeddings, it suffices to show that for any
  $\mathcal{S}$, (\ref{14}) has a unique solution in
  $\mathfrak{B}(\nu,n,\mathcal{S})$.  The proof of this property is
  very close to that of Lemma~\ref{L12}, after adaptations of the
  inequalities to the new norms, which are explained in the Appendix,
  \S\ref{Proof37}.
\end{proof}

\subsection{Ramification of  $\mathbf{F}$ at $\bf p=0$ and global properties} \label{703}
 We define $\mathfrak{B}(\nu,n,\epsilon_1)$ to be the Banach space of functions
  defined on $S_{\epsilon_1}^d=\{\mathbf{p}:\max_{i}|p_i|\le \epsilon_1\}$
  in the norm (\ref{normnu,n}) with $\mathcal{S}$ replaced by
  $S_{\epsilon_1}^d$.
\begin{Lemma}\label{Deco}
 Let 
\begin{equation}
    \label{ramF}
    {G}(\mathbf{p})=\sum_{\mb 0\preceq \mb j\prec \mb N}
p_1^{\frac{j_1}{N_1}}\cdots p_{{d}}^{\frac{j_{{d}}}{N_{{d}}}}{A}_{j_1,...,j_{{d}}}(\mathbf{p})
  \end{equation}
  where $ {A}_{j_1,...,j_{{d}}}$ are analytic at $\mathbf{p}=0$. Then the
  functions $ {A}_{j_1,...,j_{{d}}}$ are unique and for some constants
  $C_1$ and $C_2$ and large $\mathbf{p}$ we have
  \begin{equation}
    \label{est45}
    |{A}_{j_1,...,j_{{d}}}(\mathbf{p})|\le C_1|\mathbf{p}|^{C_2}\max_{\mb 0\preceq \mb j\prec\mb N} \left|G(p_1e^{2j_1\pi i},...,p_de^{2j_d\pi i})\right|
  \end{equation}
 In particular, in $S_1^d$ we have, for some constants
  $C_3$ and $C_4$,
\begin{multline}
    \label{est47}
C_3
\max_{\mb 0\preceq \mb j\prec\mb N} \sup_{|\mathbf{p}|\in S_1^d}\left|G(p_1e^{2j_1\pi i},...,p_de^{2j_d\pi i})\right|\\
  \le  \sup_{|\mathbf{p}|\in S_1^d}|{A}_{j_1,...,j_{{d}}}(\mathbf{p})|\le C_4
\max_{\mb 0\preceq \mb j\prec\mb N} \sup_{|\mathbf{p}|\in S_1^d}\left|G(p_1e^{2j_1\pi i},...,p_de^{2j_d\pi i})\right|
  \end{multline}
\end{Lemma}
\begin{Remark}\label{R00}
    We note that in (\ref{est45}) the order of analytic continuations is
  immaterial.
\end{Remark}
\begin{proof}
  The proof is by induction on ${{d}}$.  We take ${{d}}\ge 1$, assume
  (\ref{Deco}) with $A_{\mb j}$ analytic and write
  $\mathbf{p}=(p_1,\mb p^\perp)$. We have
\begin{multline}
    \label{ramF1}
    {G}(\mathbf{p})=\sum_{0\le j_{1}<N_{1}} p_{1}^{\frac{j_1}{N_1}}\left(\sum_{\{j_m<N_m;m=2,...,d\}}
p_2^{\frac{j_2}{N_2}}\cdots p_{d}^{\frac{j_{d}}{N_{d}}}{A}_{j_1,...,j_{{d}}}(\mathbf{p})\right)\\=:\sum_{0\le j_{{1}}<N_{{1}}} p_{{1}}^{\frac{j_{{1}}}{N_{{1}}}}G_{j_{{1}}}(p_1,\mb p^\perp)
  \end{multline}
  (with the convention that $G_{j_1}=A_{j_1}$ if $d=1$). We  write the system
\begin{equation}
  \label{rot}
  {G}(p_1e^{2 k\pi i},\mb p^\perp)=\sum_{0\le j_{{1}}<N_{{1}}}e^{2 kj_{{1}}\pi i/N_{{1}}} p_{{1}}^{\frac{j_{{1}}}{N_{{1}}}}G_{j_{{1}}}(p_1,\mb p^\perp);\ \ 
k=0,1,...,N_{{1}}-1
\end{equation}
which has nonzero Vandermonde determinant, from which
$G_{j_{{1}}}(p_1,\mb p^\perp)$ are uniquely determined, which in turn,
by the induction hypothesis determine $ {A}_{j_1,...,j_{{d}}}$, with
the required estimates.
\end{proof}

\begin{Lemma}\label{dec1}
  Under the assumption 1  and condition 1, the solution in
Lemma \ref{L12} can be decomposed as follows:   
  \begin{equation}
    \label{ramF2}
    \mathbf{F}(\mathbf{p},t)=\sum_{\mb 0\preceq \mathbf{j}\prec \mathbf{N}}
p_1^{\frac{j_1}{N_1}}\cdots p_d^{\frac{j_d}{N_d}}\mathbf{A}_\mathbf{j}(\mathbf{p},t)
  \end{equation}
  where $
  \mathbf{A}_\mathbf{j}(\mathbf{p}, t)\in\mathfrak{B}(\nu,n,\mathcal{S})$
  are analytic at $\mathbf{p}=0$. Furthermore, in analyzing the
  continuations in restricted sectors $\mb{p}e^{2\pi i
    \mathbf{j}}\in\mathcal{S}_\phi$ we have for some $\nu$, in the
  norm defined in (\ref{16} ) (cf. also Remark~\ref{R00})
  \begin{equation}
    \label{trueexp}
   \max\left\{\|\mathbf{F}(\cdot e^{2\pi i \mathbf{j}},\cdot)\|_\nu,
    \{\|\mathbf{A_j}(\cdot,\cdot) \|_{\nu}\};\mb 0\preceq \mathbf{j}\prec 
\mathbf{N}\right\}=K<\infty
  \end{equation}
\end{Lemma}
\begin{proof}
  We consider the equation (\ref{14}) on $\mathfrak{B}(\nu,n,
  \mathcal{S})^{\tilde{N}}$ where $\tilde{N}$ counts the
  $\mathbf{A}_{\bf j} (\cdot,t)$ via the decomposition
  (\ref{ramF2}).  Noting that
\begin{equation}\label{conviden}
  p^{\alpha}*p^{\beta}=
\frac{\Gamma(\alpha+1)\Gamma(\beta+1)}
{\Gamma(\alpha+\beta+2)}p^{\alpha+\beta+1}
\end{equation}
it is straightforward to show that the space of functions of the form
(\ref{ramF}) is stable under convolution. Since $R({\bf p}, t)$ and
therefore ${\bf F}_0 ({\bf p}, t)$ are of the form (\ref{ramF2}) it
follows that $\mathcal{N}$ leaves the space of $\mathbf{F}$ of the
form (\ref{ramF2}) invariant. Using the estimates (\ref{est47}) we see
that $\mathcal{N}$ is well defined in a small ball of radius
$\epsilon_2$ in in $\mathfrak{B}(\nu,n,\mathcal{S})$ and that it is a
contraction there.  Therefore the solution to (\ref{14}) is of the
form (\ref{ramF2}). For $\mb{p}e^{2\pi i \mathbf{j}}\in\mathcal{S}_\phi$,
$\|\mathbf{F}(\mathbf{p}e^{2\pi i \mathbf{j}})\|_\nu$ are well
defined.  Using again Lemma~\ref{Deco} the first statement follows. To
show finiteness of $\|\mathbf{A}_{\bf j} (\cdot,t) \|_{\nu} $ it
suffices to prove finiteness of $\|\mathbf{F}(\mathbf{p}e^{2\pi i
  \mathbf{j}})\|_\nu$. To this end, we note that all these
continuations satisfy equations of the type (\ref{eq:defconv}) with
coefficients satisfying the requirements in \S\ref{S3} and thus the
result follows from Lemma~\ref{L12}.

\end{proof}
\begin{Lemma}\label{ident2}
  Assume $\mb G$ is an entire function of exponential order $n$, more
  precisely satisfying the inequality $|\mb G(\mb p)|\le Ce^{\nu|\mb
    p|^n}$ for some constants $C,\nu$ and that in a sector
  $\mathcal{S}_\phi=\{\mb p:|\mb p|>0,\max_{i}|\arg(p_i)|<\phi\}$, it
  grows at most exponentially, $|\mb G(\mb p)|\le Ce^{\nu_1|\mb p|}$.
  Then there exists a function $\mb G_1$ increasing at most
  exponentially $|\mb G_1(\mb p)|\le Ce^{\nu_2|\mb p|}$ in any proper
  subsector of $\mathcal{S}_{\phi_1}$ where $\phi_1 = \frac{\pi}{2
    (n-1)} + \frac{n \phi}{n-1}$ and such that $\mb G(\mb z^n)$ is
  analytic at $\mb z=0$, such that
  \begin{equation}
    \label{ident}
    \mb g(\mb x):=\int_0^{\infty}e^{-\mb p\cdot \mb x}\mb G(\mb p)d\mb p=\int_0^{\infty}e^{-\mb p\cdot \mb x^{\frac{n}{n-1}}}\mb G_1(\mb p)d\mb p
  \end{equation}
\end{Lemma}
\begin{proof}
  We start with the case when $\mb G$, $\mb x$ and $\mb p$ are scalar,
  the general case following in a quite straightforward way as
  outlined at the end.
  
  The assumptions on $G$ ensure that the first integral in
  (\ref{ident}) exists and $g(x)$ has an asymptotic power series in
  powers of $x^{-1}$ in a sector of opening $\pi+2\phi$ centered on
  $\RR^+$. The function $g_1(x)=g(x^{(n-1)/n})$ has a (noninteger)
  power series asymptotics in a sector of opening
  $\frac{n}{n-1}(\pi+2\phi)$ and by the general theory of Laplace
  transforms, $G_1:=\mathcal{L}^{-1}g_1$ is analytic in a sector of
  opening $\frac{n}{n-1}(\pi+2\phi)-\pi$ centered on $\RR^+$, Laplace
  transformable, with Laplace transform $g_1$. It follows that
\begin{equation}
  \label{formulaG1}
  G_1(p)=
\frac{1}{2\pi i}\int_{c-i\infty}^{c+i\infty}
e^{pu}
\int_0^{\infty}e^{-qu^{(n-1)/n}}G(q)dqdu=:\int_0^{\infty}K_{\frac{n-1}{n}}(p,q)G(q)dq
\end{equation}
We show that $G_1$ has a convergent expansion in powers of $p^{1/n}$
at zero. The function
\begin{equation}
  \label{resK}
  K_{\frac{n-1}{n}}(p,q)=\left(\frac{q}{p}\right)^n
C_{\frac{n-1}{n}}(q^n/p^{n-1})
\end{equation}
is \'Ecalle's acceleration kernel
\cite{Balser,[E6]}.  For $\alpha\in(0,1)$, with $\beta=1-\alpha$,
$c=\beta\alpha^{\alpha/\beta}$, the function $C_\alpha$ is an entire
function and has the following asymptotic behavior \cite{Balser,[E6]}:
\begin{equation}
  \label{Ecalleacc2}
 C_{\alpha}(x)\sim \frac{\alpha^{\frac{1}{2\beta}}}{\sqrt{2\pi\beta}} x^{1/2}e^{-cx};\ \ |x|\to\infty, \ |\arg x|<\frac{\pi}{2}
\end{equation}
Using (\ref{resK}) we see that
\begin{equation}
  \label{resK2}
  \int_0^{\infty}K_{\frac{n-1}{n}}(p,q) q^k dq =p^{(nk-k-1)/n}\int_0^\infty s^{k+n} C_{\frac{n-1}{n}} (s^n) ds
\end{equation}
We expand the entire function $G$ in series about the origin, 
$G(q)=\sum_{k=1}^{N-1}g_k q^k+R_N(q)$
and note that
\begin{equation}
  \label{eq:RN}
  |R_N(q)|\le \sum_{k=N}^{\infty}|G^{(k)}(0)||q|^k/k!\le \sum_{k=0}^{\infty}|G^{(k)}(0)||q|^k/k!\le
Ce^{\nu_5 |q|^n}=E(q)
\end{equation}
uniformly in $\CC$. By (\ref{Ecalleacc2}) and (\ref{eq:RN}) 
 $E(q) C_{\alpha}(q^n/p^{n-1})$ is, for small enough $p$,  in $L_1[0,\infty]$ 
in $q$.
By dominated convergence, we have
$$\int_0^{\infty}K_{\frac{n-1}{n}}(p,q)G(q)dq=\lim_{N\to\infty}\int_0^{\infty}K_{\frac{n-1}{n}}(p,q)\sum_{k=1}^{N-1}g_k
q^k dq$$
and, using (\ref{resK2}) it follows that for small $p$, $G_1$
is the sum of a convergent series in powers of $p^{1/n}$, as
stated\footnote{To estimate the radius of convergence of this series
  it is  convenient to start from the duality (\ref{ident}) and
  apply Watson's lemma, using  Cauchy's formula on a
  circle of radius $k^{1/n}/(n\nu)^{1/n}$ to bound $|\mb G^{(k)}(0)|$ .}.

The argument for $d$ variables and vectorial $\mathbf{G}$ is nearly
the same: a vectorial $\bf G$ is treated componentwise, while the
assumptions ensure that the multidimensional integrals involved can be
taken iteratively, the estimates being preserved in the process.

\end{proof}

Collecting the results of Lemma~\ref{dec1} and Lemma~\ref{ident2}
applied to each of the $\bf A_j$, the proof of Theorem~\ref{TrB}
follows.

\begin{Note}\label{Opti2}
  In the example $\partial_tu+(-\partial_x)^nu=0$ we have
  $\phi=\frac{\pi}{2n}$. Formal exponential solutions have the
  behavior, to leading order, $\exp\left(c_n
    (-x)^{\frac{n}{n-1}}t^{-\frac{1}{n}}\right)$ with
  $c_n=(n-1)/4/n^{\frac{n}{n-1}}$ (for all determinations of
  $(-x)^{\frac{n}{n-1}}$).  This also points to $x^\frac{n}{n-1}$ as
  natural variable and indicates that the sector of summability cannot
  be improved since it is bordered by (anti)stokes lines.
 
\end{Note}

\section{Short time existence and asymptotics in special cases}
\label{smalltime}

In some cases, the Borel summation approach can be adapted to study
short time existence of sectorial solutions and study small time
asymptotics.  One important application is in the analysis of
singularity formation in PDEs \cite{CPAM3}.  For simplicity,
and since some assumptions are less general than in the rest of the
paper, we restrict to $d =1$ (scalar case) in this section. 

 We
motivate the assumptions made by looking at a particular example
arising in Hele-Shaw flow with surface tension
\begin{equation} 
\label{eqmhdym} 
H_t= -\frac{H^3}{2} + H^3 H_{zzz},\ \   H(z, 0) = z^{-1/2} 
\end{equation}
the modified  Harry-Dym equation (see \cite{TNV}, where it
arises with $\xi=z+t$ (as  a local approximation near an initial zero
of the derivative of a conformal mapping).  
\subsection{Formal series, preparation of normal form.}\label{FSP}

\z {\bf Note:} To simplify notation, in the following we let
$\mathfrak{p}$ stand for generic {\em polynomials}, $\mathfrak{p}^+$
for polynomials with {\em nonnegative coefficients}, and
$\mathfrak{p}_{(n)}$ for polynomials of {\em degree $n$}. Similar
conventions are followed for $\mathfrak{h}$ which represents {\em
  homogeneous} polynomials.  Substituting in (\ref{eqmhdym}) a
power-series of the form $\sum_{n=0}^\infty t^n H_n (z) $
where $H_0 = z^{-1/2}$ yields the recurrence
\begin{equation}
\label{recurrence}
(n+1) H_n = - \frac{1}{2} \sum_{n_j \ge 0, \sum_{j=1}^3 n_j = n} \!\!\!
H_{n_1} H_{n_2} H_{n_3} 
+\!\!\! \sum_{n_j \ge 0, \sum_{j=1}^4 n_j = n} \!\!\!H_{n_1} H_{n_2} H_{n_3} 
H_{n_4}^{\prime \prime \prime} 
\end{equation}
which  inductively shows that $H_n= z^{-1/2} \mathfrak{h}_{(n)} (
z^{-9/2}, z^{-1} )$. We let
\begin{equation}
\label{A0.2.3}
g_N(x,t):=\sum_{k=0}^Nt^nH_n(z)=x^{-1/3}\sum_{n=0}^N \mathfrak{h}_{(n)}(tx^{-3},tx^{-2/3});\text{ where }x = \frac{2}{3} z^{3/2}
\end{equation}
In terms of $x$, (\ref{eqmhdym}) becomes,
\begin{equation}
\mathcal{N} (H):=H_t+\frac{1}{2}H^3-\frac{3x}{2}H^3 H_{xxx}
-\frac{3}{2}H^3H_{xx}+\frac{1}{6x}H^3H_x=0
\end{equation}
 It is straightforwardly shown that
\begin{equation}
\label{eq:eqA02}
\mathcal{N} g_N(x,t)
= t^{-1} x^{-\frac{1}{3}} \mathfrak{p}_{(4 N+1)} (t x^{-3}, t x^{-2/3} )
\end{equation}
where for small $x_1,x_2$ we have moreover
\begin{equation}
  \label{hm2}
\mathfrak{p}_{(4 N+1)} (x_1,x_2) = \mathfrak{h}_{(N+1)} (x_1,x_2) \left [ 1 + O(x_1,x_2)  \right ]
\end{equation}
It is then natural to substitute :
\begin{equation}
\label{eqHf}
~H(z(x) ,t) = g_N(x, t) + x^{-2} f(x,t) 
\end{equation} 
 into (\ref{eqmhdym}); we choose  without  loss of generality $N \ge 3$.

 It will
 follow from the analysis that $ |f(x, t)| =o\left ( x^{5/3}
   \mathfrak{h}_{(N)} \left ( t x^{-3}, t x^{-2/3} \right )\right ) $
 for small $t^{1/3}x^{-1}$ with $\arg x \in \left (-\frac{\pi}{2} -
   \phi, \frac{\pi}{2} + \phi \right ) $ and
 $\phi\in(0,\frac{\pi}{6})$, thus $H\sim  \sum_{n=0}^\infty t^n H_n (z) $ 
for small $t^{1/3}x^{-1}$ (see Corollary~\ref{mhdym}).
 
 Substitution shows that $f (x, t)$ satisfies an equation of the form
 (\ref{1.a}), with $n=3$ (third order, $m=1$ (scalar case), with (cf.
 also (\ref{3}), and (\ref{valpha}) below)
\begin{equation}
\label{eq:eqA2}
r(x, t) = t^{-1} x^{5/3} \mathfrak{p}_{(4 N+1)} (t x^{-3}, t x^{-2/3});\ \ 
b_{{\bf q}, k} = x^{-\beta k} \sum_{j=1}^{J_{\bf q}} 
x^{-\alpha_{{\bf q},k}} 
\mathfrak{p}_{{\bf q}, k; j} (t x^{-3}, t x^{-2/3} ) 
\end{equation} 
{\bf Note:} By (\ref{hm2}), $r(x, t)$ is small for small $t$ or large
$x$, in spite of the prefactor $t^{-1} x^{5/3}$.
 
\subsection{More general setting.}\label{MGS}

\z {\bf Setting 1.}
We take $\rho_0 = 0$, suitable for algebraic initial
conditions in the domain, and consider the domain $\mathcal{D}_{\phi,
  0, x} $, with $\phi < \frac{\pi}{2n}$ small enough to ensure
(\ref{condnew}). Taking ${\bf f}
(x, t) - {\bf f}_I (x) $ as the unknown function we may assume
$$ {\bf f}_I (x) = {\bf 0} $$
(see Note 3 after Theorem~\ref{Tasympt}) and require that
\begin{equation}
\label{eq:eqr}
| {\bf r} (x, t) | \le   
t^{-1}\sum_{j=1}^{J_r} |x|^{\omega_{j}} \mathfrak{h}^+_{(n'_j)} \left ( 
t^{\gamma_1} |x|^{-\beta_1}, ....,t^{\gamma_K} |x|^{-\beta_K} \right )
\end{equation}
where the degrees $n'_j$ satisfy 
\begin{equation}
\label{eq:eqn'}
{n'}_j \beta_l - \omega_j \ge  1 ,\\ ~{\rm for} ~1 \le l \le K, ~1 \le j \le J_r
\end{equation}
 (As before, (\ref{eq:eqn'}) implies
that $r(x, t)$ is small for large $x$ or small $t$).
The positive constants $\omega_1$, $\omega_2$, ...,$\omega_{J_r}$,
$\beta_1$, $\beta_2$, ...,$\beta_K$ and $\gamma_1$, $\gamma_2$,
...,$\gamma_K$,  are  restricted by the
condition
\begin{equation}
\label{eqhatn}
{\hat n} :=\frac{\beta_1}{\gamma_1} \ge n 
\end{equation}
The labeling is chosen so that
\begin{equation}
\label{eqorder}
{\hat n} =\frac{\beta_1}{\gamma_1} \ge \frac{\beta_2}{\gamma_2} 
....\ge \frac{\beta_K}{\gamma_K} 
\end{equation}
Also, if for some $1 \le j \le K-1$, $
\frac{\beta_j}{\gamma_j} = \frac{\beta_{j+1}}{\gamma_{j+1}} $, we
arrange $\beta_{j} > \beta_{j+1}$.  The $\omega_j$ are arranged
increasingly:
\begin{equation}
\label{eqomorder}
\omega_1 < \omega_2 < .... < \omega_{J_r} 
\end{equation}
Furthermore, 
for any $x \in \mathcal{D}_{\phi,0, x}$, we require
\begin{equation}
\label{eq:eqbqk}
| {\bf b}_{{\bf q}, {\bf k}} (x, t) | \le  |x|^{-\beta | {\bf k} | } 
\sum_{j=1}^{J_{\bf q}} |x|^{-\alpha_{{\bf q}, j}} \mathfrak{p}^+_{{\bf q}, {\bf k}, j} 
\left ( 
t^{\gamma_1} |x|^{-\beta_1}, ....,t^{\gamma_K} |x|^{-\beta_K} \right )
\end{equation}
\begin{equation}
\label{eq:alphaqorder}
\beta > 0 ~,~\\ ~~~\alpha_{{\bf q}, 1} > \alpha_{{\bf q}, 2} > ... > 
\alpha_{{\bf q}, J_{\bf q}} ~~;\ ~~
\mb b_{\bf q, k}\ne 0\Rightarrow \alpha_{\bf q, j} + \beta |{\bf k} | \ge 0 
\end{equation}
If only finitely many 
${\bf b}_{{\bf q}, {\bf k}}$ are nonzero we allow
\begin{equation}
\label{betage0}
\beta \ge 0  
\end{equation}
We also require that for all ${\bf q}$, ${\bf k}$ for which
${\bf b}_{\bf q, k} \ne 0$ we have
\begin{equation}
\label{eqcond}
m_{\bf q,k}:={\hat n} + 
\omega_1 ( |\mb q| - 1)
-\alpha_{{\bf q}, 1} + (\omega_1 - \beta) | {\bf k} | - 
\frac{\hat n}{n} \sum_{j,l} j q_{l,j} \ge 0
\end{equation}

{\bf Note:} Assumption (\ref{eqcond}) is satisfied by modified
Harry-Dym and by certain classes of nonlinear PDEs and initial
conditions-- for instance, the thin-film equation $h_t + (h^3
h_{xxx})_x = 0$, with singular initial condition $h(x, 0) =
x^{-\alpha}$ for $\alpha > 0$, but is generally quite restrictive.
Weakening it requires more substantial modifications of the framework
and will not be discussed here.

{\bf Setting 2.}  Better properties are obtained under the assumptions
described below.
\bigskip
\begin{equation}
  \label{efb}
  \begin{tabular}{lllllll}
$ {\hat n} = n $\\ \\
$ \mathcal{P} (-s) = s^n $\\ \\
$ {\bf r} (x, t)  = 
\frac{1}{t} \sum_{j=1}^{J_r} x^{\omega_{j}} {\bf \mathfrak{a}}_{j} \left ( 
t^{\gamma_1} x^{-\beta_1}, ....,t^{\gamma_K} x^{-\beta_K} \right ) $\\ \\
$ {\bf b}_{{\bf q}, {\bf k}} (x, t)  = x^{-\beta  |{\bf k}|  } 
\sum_{j=1}^{J_{\bf q}} x^{-\alpha_{{\bf q}, j}} {\bf \mathfrak{a}}_{{\bf q}, {\bf k}, j} 
\left ( t^{\gamma_1} x^{-\beta_1}, ....,t^{\gamma_K} x^{-\beta_K} \right ) $
\end{tabular}
\end{equation}

\bigskip

\z where ${\bf \mathfrak{a}}_j$, ${\bf \mathfrak{a}}_{{\bf q}, {\bf k},
  j}$ are {\em analytic} near the origin and for
small $|\bf z|$ we require, with the same restriction (\ref{eq:eqn'}) on $n'_j$,
\begin{equation}
\label{eq:eqQnew}
|\mathfrak{a}_{j} (\mb z)| \le \mathfrak{h}^+_{(n'_j)}(|z_1|,...,|z_n|)
\end{equation}
The restrictions on the numbers $\beta_1$, $\beta_2$, ...$\beta_K$,
$\gamma_1$, $\gamma_2$, ...$\gamma_K$, $\alpha_{{\bf q}, j}$, etc.
are as in Setting 1. Furthermore, we assume that there is an
$\omega\in\RR^+$ so that the nonnegative numbers
\begin{equation}\label{eqrational}
m_{\bf q,k},\ \omega_2 - \omega_1,...,\omega_{J_r} - \omega_1,\ \alpha_{{\bf q},1}- \alpha_{{\bf q}, 2},..., \alpha_{{\bf q}, 1}-\alpha_{{\bf q}, J_{\bf q}},\ n
\gamma_2-\beta_2,...,n\gamma_K - \beta_K
\end{equation}
are {\em integer multiples of $n\omega$}. This condition, satisfied for the
problem (\ref{eqmhdym}), comes out naturally in a number of examples
and ensures the existence of a ramified variable in which the
solutions are analytic. We choose $\omega>0$ to be the largest with
the property above. Define
\begin{equation}
\label{eq:eqzeta}
\zeta = y t^{-1/{n}} ~,~~ 
\\ {\hat {\bf f}} (\zeta, t) = {\bf f} (t^{1/{n}} \zeta, t)
\end{equation}
and
\begin{equation}
\label{eqhatD}
{\hat D}_{\phi, \rho} = \left \{ \zeta: |\zeta| > \rho;\ |\arg \zeta|<\phi \right \} 
\end{equation} 
\begin{Theorem}
\label{Tasympt} (i) In Setting 1, under Assumption 1, 
there exists for large enough
$\rho$ a unique solution ${\bf {\hat f}} (x t^{-1/{\hat n}}, t)$ 
to (\ref{1.a}),
for
$\zeta = x t^{-1/{\hat n}} \in {\cal {\hat D}}_{\phi, \rho}$ 
and, with $n'_j$ as in (\ref{eq:eqn'}),
\begin{equation}
\label{fsbound}
|{\bf {\hat f}} (\zeta, t)| \le \sum_{j=1}^{J_r} |\zeta|^{\omega_j} t^{\omega_j/{\hat n}}
\mathfrak{h}_{(n'_j)} (|\zeta|^{-\beta_1}, t^{\gamma_2-\beta_2/{\hat n}} |\zeta|^{-\beta_2}, ...,
t^{\gamma_K-\beta_K/{\hat n}} |\zeta|^{-\beta_K} ) 
\end{equation}
(ii) In Setting 2, under Assumption 1, for any $T>0$ there is a
$\rho=\rho(T)>0$ so that the mapping
$$(\zeta,\theta)\to \theta^{ -\frac{\scriptstyle \omega_1}{\scriptstyle n\omega}}
\hat {\mb f}(\zeta,\theta^{1/\omega})$$
is analytic in $\mathcal{\hat
D}_{\phi, \rho}\times \{\theta:|\theta|<T\}$.
\end{Theorem}

\noindent{\bf Notes:} 1. The function $\rho$ will, generally, increase with $T$.

2. The restriction $d=1$ is not essential, but made for the sake of
simplicity.

3. In these settings, there is a duality between large $x$ and small
$t$ in the asymptotics: $\zeta$ can be large either due to largeness
of $x$  or smallness of $t$.  For $t $ in a fixed interval,
there exists some $\rho$ so that the asymptotic bounds are satisfied
for $\zeta \in {\hat D}_{\phi,\rho}$.

4. The following example shows that the requirement ${\hat n} \ge n $
is natural. In the equation $g_t +(-\partial_x)^n g= 0$ with $g(x,0) =
x^{-\alpha}$, substituting the expansion $g(x, t) = x^{-\alpha} +
\sum_{n\in\NN}t^ng_n(x)$, we get $g_n (x) = O(x^{-\alpha -n})$. Thus
one of the scales that emerge in the formal expansion is $t/x^n$.  On
the other hand, in view of (\ref{eq:eqr}) and (\ref{eq:eqbqk}) the
most singular term as $x \rightarrow 0$  is of the order $t/x^{\hat
  n}$ since ${\hat n} =\frac{\beta_1}{\gamma_1}$.  Combining with the
above discussion we see that $\hat n \ge n$.

5.  The leading order term in the  Taylor expansion of
$\theta^{-\frac{\omega_1}{n\omega}} {\bf \hat f}$, ${\bf {\hat f}_0}$,  satisfies
an easily obtained ODE.  The convergence of the series in part (ii)
implies that singularities of ${\bf {\hat f}_0}$ can be related to
actual singularities of the PDE for small time and this is the subject
of another paper (\cite{CPAM3}).

\begin{Corollary}
\label{mhdym} For  the initial value problem (\ref{eqmhdym}),
 for any $T>0$ there is a $\rho=\rho(T)$ such that 
\begin{equation}
  \label{asHD}
  H(z,t)=\sum_{k=0}^{\infty}t^{\frac{7k+1}{9}}G_k(z t^{-2/9})
\end{equation}
where the series converges in the region $\{(z,t):|t|<T,|z|>\rho,|\arg
z|<\frac{4}{9} \pi \}$ and $G_k (\zeta)$ 
are analytic in the sector $\{\zeta:|\zeta|>\rho,
|\arg \zeta|<\frac{4}{9} \pi \}$.
\end{Corollary}

\subsection{Proof of Theorem \ref{Tasympt} (i)}

It is convenient to make rescalings of variables in Borel space as well. We note that 
\begin{equation}
\label{eqlaplace}
{\hat {\bf f}} (\zeta, t) = t^{-1/{\hat n}}
\int_0^\infty e^{-s \zeta} {\hat {\bf F}} (s, 1; t)
ds
\end{equation}
where 
\begin{equation}
\label{eq:eqzeta2}
s = p t^{1/{\hat n}} ~,~ \\ {\hat {\bf F}} (s, \lambda; t) =
{\bf F} (t^{-1/{\hat n}} s, t \lambda)
\end{equation}
We use similar rescaling to define ${\hat {\bf R}} (s, \lambda; t), {\hat {\bf B}}_{{\bf q}, {\bf k}} (s, \lambda; t)$ and ${\bf {\hat F}_0} (s,\lambda; t)$
where now
\begin{equation}
\label{eq:eqhatF0}
{\bf {\hat F}_0} (s, \lambda; t) = t \lambda \int_0^1 e^{-t \lambda \mathcal{P} (-s t^{-1/{\hat n}})
 (1-\tau)}
{\bf {\hat R}} (s, \lambda \tau; t) d \tau
\end{equation}
We let $\mu_{\bf q,k}=1-\hat{n}^{-1} \left(| {\bf q} | + | {\bf k} |+\sum_{j=1}^n\sum_{l=1}^m j
q_{l,j}\right)$. Using (\ref{14}),
straightforward calculations show that
\begin{multline}\label{eq:eqhatF}
{\bf {\hat F}} (s,\lambda; t)={\hat {\mathcal N}} ({\bf {\hat F}}) (s, \lambda; t)
\equiv
{\bf {\hat F}_0} (s, \lambda; t) +
{\sum_{{\bf q}
\succeq 0}}' \sum_{{\bf k} \succeq 0}
\lambda t^{\mu_{\bf q,k} } \\ \times 
\int_0^1
e^{-t \lambda \mathcal{P} (-s t^{-1/{\hat n}})
 (1-\tau)}
\left \{ {\bf {\hat B}}_{{\bf q},{\bf k}} * {\bf \hat{F}}^{*{\bf k}}
* \sideset{^*}{}\prod_{l=1}^{m} \sideset{^*}{}\prod_{j=1}^{n}
\left ( (-s)^j {\hat F}_l \right )^{*q_{l,j}} \right \}
(s, \lambda \tau, t)
d\tau
\end{multline}
With slight abuse of notation we drop the hats from the newly defined
functions. Let now
\begin{equation}\label{eqhatD1}
{{\cal S}}_{\phi} \equiv
\left \{ s : \arg s \in
  (-\phi, \phi ),~
0 < |s| < \infty ,~
 0 < \phi < \frac{\pi}{2n} \right \}
\end{equation}
and consider the Banach space  ${ {\cal A}}_\phi $ of analytic functions in
${{\cal S}}_\phi $, continuous in $\overline{S}_\phi$ in the norm
\begin{equation}\label{normhat}
\| { {\bf F}} (\cdot,\cdot; t) \|_{\nu} = \sup_{0 \le \lambda \le 1, s \in
{{\cal S}}_\phi} (1+|s|^2) e^{-\nu |s|} | { {\bf F}} (s, \lambda; t) |
\end{equation}
\begin{Lemma}\label{LL37}
With ${\bf r} (x, t)$ satisfying (\ref{eq:eqr}) we have 
$$\| {\bf { F}}_0 (., .; t) \|_\nu 
\le 
e^{a t} \sum_{j=1}^{J_r} \nu^{\omega_j +1} t^{(\omega_j+1)/{\hat n}} 
\mathfrak{h}^+_{{n'}_j} 
\left ( \nu^{-\beta_1}, 
t^{\gamma_2-\beta_2/{\hat n}} \nu^{-\beta_2},..., 
t^{\gamma_{K} 
- \beta_{K}/{\hat n}} \nu^{-\beta_K} \right )$$
for $\nu$ large (independent of $t$ for small $t$), where $-a$ is
the lower bound of $\Re \mathcal{P} (p)$. 
\end{Lemma}

\begin{proof}
>From (\ref{eq:eqr}), (\ref{eq:eqn'}) and 
applying Lemma \ref{L1} (with $\rho=0$; see Remark \ref{rem1}) we have
\begin{multline*}
   | {\bf { R}} (s, \lambda; t) |\\ \le 
\frac{1}{t \lambda} \sum_{j=1}^{J_r} |s|^{-\omega_j -1} t^{(\omega_j+1)/{\hat n}}
\mathfrak{h}^+_{n'_j} 
\left (\lambda^{\gamma_1} |s|^{\beta_1}, 
\lambda^{\gamma_2} t^{\gamma_2-\beta_2/{\hat n}} |s|^{\beta_2},..., 
\lambda^{\gamma_K} t^{\gamma_{K} 
- \beta_{K}/{\hat n}} |s|^{\beta_K} \right ) 
\end{multline*}
For $\lambda\in (0,1)$ we have 
$\left | e^{-t \mathcal{P} (-s t^{-1/{\hat n}} ) \lambda (1-\tau)} 
\right | \le e^{a t}$ and thus (cf. (\ref{eq:eqhatF0}))
\begin{multline}
\label{eq:eqF0b}
| {\bf { F}_0 } (s, \lambda; t) | \\\le e^{a t}
\sum_{j=1}^{J_r} |s|^{-\omega_j -1} t^{(\omega_j+1)/{\hat n}}
\mathfrak{h}^+_{n'_j} 
\left (\lambda^{\gamma_1} |s|^{\beta_1}, 
\lambda^{\gamma_2} t^{\gamma_2-\beta_2/{\hat n}} |s|^{\beta_2},..., 
\lambda^{\gamma_K} t^{\gamma_{K} 
- \beta_{K}/{\hat n}} |s|^{\beta_K} \right )
\end{multline}
Bounding each term
of the polynomial ${\mathfrak h}^+_{n'_j}$ in $\|\cdot\|_\nu$ we obtain  
$$
\| {\bf {\hat F}_0} (., . ; t) \|_\nu \le e^{a t} \sum_{j=1}^{J_r}
\nu^{\omega_j +1} t^{(1+\omega_j)/{\hat n}} {\mathfrak h}^+_{n'_j}
\left ( \nu^{-\beta_1}, t^{\gamma_2-\beta_2/{\hat n}}
  \nu^{-\beta_2},..., t^{\gamma_{K} - \beta_{K}/{\hat n}}
  \nu^{-\beta_K} \right )$$
The proof now follows, choosing $\nu$
sufficiently large and using (\ref{eq:eqn'}) and (\ref{eqorder}),
(\ref{eqomorder}).
\end{proof}

\begin{Lemma}
\label{L30} For large $\nu$, we have
$$ \| {\bf { B}}_{{\bf q}, {\bf k}} * {\bf { F}} \|_\nu
\le    
{ c}_{{\bf q}, {\bf k}} (\nu, t)  \| {\bf { F}} \|_\nu  ,{\rm where}
$$
\begin{equation}
  \label{defck}
c_{\bf 0, 0} = \mb 0;\ \ { c}_{{\bf q}, {\bf k}} (\nu, t)  =  
\nu^{-\beta | {\bf k} | } 
t^{(1-\beta |{\bf k}| )/{\hat n}} 
\sum_{j=1}^{J_{\bf q}} K_j \nu^{-\alpha_{{\bf q},j}} 
t^{-\alpha_{{\bf q}, j}/{\hat n}}\ ((\bf q,k)\ne 0)
\end{equation} 
with $K_j$ constants independent of ${\bf q}$, ${\bf k}$, $\nu$ and
$t$. 
\end{Lemma}
\begin{proof}
  Note first that ${\bf b}_{{\bf 0}, {\bf 0}} = \bf 0$ hence $c_{\bf
    0, 0} = \bf 0$. From (\ref{eq:eqbqk}) and Lemma \ref{L1} (with
  $\rho=0$),
$$ \left | {\bf B}_{{\bf q}, {\bf k}} (p, t)\right |  
\le |p|^{\beta |{\bf k}| - 1} \sum_{j=1}^{J_{\bf q}} 
|p|^{\alpha_{{\bf q},j}} \mathfrak{p}^+_{{\bf q}, {\bf k}, j} 
\left ( t^{\gamma_1} |p|^{\beta_1},  
t^{\gamma_2} |p|^{\beta_2},...,t^{\gamma_K} |p|^{\beta_K} \right ) $$  
Switching from $(p, t)$ to $(s, \lambda; t)$, 
\begin{multline*}
| {\bf { B}}_{{\bf q}, {\bf k}} (s, \lambda; t) | \le t^{(1-\beta |{\bf k}| ) /{\hat n}} 
|s|^{\beta | {\bf k} | - 1} 
\\ 
\times \sum_{j=1}^{J_{\bf q}} |s|^{\alpha_{{\bf q}, j} } 
t^{-\alpha_{{\bf q}, j}/{\hat n}}  
\mathfrak{p}^+_{{\bf q}, {\bf k}, j} 
\left (\lambda^{\gamma_1} |s|^{\beta_1}, 
\lambda^{\gamma_2} 
t^{\gamma_2-\beta_2/{\hat n}} |s|^{\beta_2},..., \lambda^{\gamma_K} t^{\gamma_{K} 
- \beta_{K}/{\hat n}} |s|^{\beta_K} \right ) 
\end{multline*}
For large $\nu $, using Lemma \ref{L4} (with $\rho = 0$) to bound in
norm the terms of $\mathfrak{p}^+_{{\bf q}, {\bf k}, j}$
\begin{multline}\label{qtilde}
\| {\bf { B}}_{{\bf q}, {\bf k}} * { {\bf F}} | \le \| {\bf  F} \|_\nu
\,\,t^{(1-\beta |{\bf k}| ) /{\hat n}}  
|\nu|^{-\beta | {\bf k} |} \\  \times
\sum_{j=1}^{J_{\bf q}} |\nu|^{-\alpha_{{\bf q}, j} } 
t^{-\alpha_{{\bf q}, j}/{\hat n}}  
\mathfrak{p}^+_{{\bf q}, {\bf k}, j} 
\left (\lambda^{\gamma_1} \nu^{-\beta_1}, 
\lambda^{\gamma_2} 
t^{\gamma_2-\beta_2/{\hat n}} \nu^{-\beta_2},..., \lambda^{\gamma_K} t^{\gamma_{K} 
- \beta_{K}/{\hat n}} \nu^{-\beta_K} \right ) 
\end{multline}
Clearly, for large $\nu$, $\mathfrak{p}^+_{{\bf q}, {\bf k}}$ can be
replaced in (\ref{qtilde}) by a constant $K_j$.  Using 
(\ref{eqorder}) and (\ref{eq:alphaqorder}) the conclusion follows.
\end{proof}

Let now
\begin{multline*}
C(\phi, T) =
\max \Bigg \{\\
\sup_{p \in \mathcal{S}_\phi, |p| > R, 0 \le l' \le n, \gamma >0} 
\left ( \frac{|p|^n}{\Re \mathcal{P}(-p)} \right )^{l'/n}
\frac{1 - e^{-\gamma}}{\gamma^{1-l'/n}}, 
\sup_{p \in \mathcal{S}_\phi, |p| \le  R,  0 \le l' \le n} 
t^{l'/n} |p|^{l'} e^{-t \Re\mathcal{P} (-p)} \Bigg\} 
\end{multline*}
where $R$ is the same as in the proof of Lemma \ref{newcond}.
\begin{Lemma}
\label{eq:eqlemmasmall}
For $\nu$ large enough, ${ {\cal N}}$ is  contractive, and thus
there exists unique solution ${\bf { F}}$ of (\ref{eq:eqhatF}).
\end{Lemma}

\begin{proof}
  For $\nu$ large enough, (\ref{eqcond}), Lemma~\ref{LL37} and
Lemma~\ref{L30} imply
\begin{equation}
\label{eq:eqfirst}
C(\phi, T) 
{\sum_{{\bf q}
\succeq 0}}' \sum_{{\bf k} \succeq 0}
t^{\mu_{\bf q,k} }
c_{{\bf q}, {\bf k}} (\nu, t)
 \| 2 {\bf { F}_0} \|^{| {\bf k} | + | {\bf q} |}
\le\| {\bf { F}_0} \|_\nu
\end{equation}
and
\begin{equation}
\label{eq:eqsecond}
C(\phi, T) 
{\sum_{{\bf q}
\succeq 0}}' \sum_{{\bf k} \succeq 0}
t^{\mu_{\bf q,k} } 
c_{{\bf q}, {\bf k}} (\nu, t) ( | {\bf q} |  +
| {\bf k} | )
\| 6 {\bf { F}_0} \|^{|{\bf k}| + | {\bf q}
| -1} \le 1
\end{equation}
 Now, Lemma \ref{L4.9.5} (with $\rho_0=0$, $d=1$ and
$s$ replacing $p$), and Lemma \ref{L30} imply
\begin{equation*}
\left | \left \{ { {\bf B}}_{{\bf q}, {\bf k} } * { {\bf F}}^{*{\bf k}}
* \sideset{^*}{}\prod_{l=1}^{m}
\sideset{^*}{}\prod_{j=1}^{n} \left (s^j { F}_l \right )^{*q_{l,j}}
\right \} (s, \lambda \tau; t) \right |
\le \frac{e^{\nu |s|} |s|^{\sum j q_{l,j} } }{
M_0 (1 + |s|^2 ) } c_{{\bf q}, {\bf k}} (\nu, t)
\| { {\bf F}} \|_\nu^{|{\bf q}| + | {\bf k} | }
\end{equation*}
Also, note that if $l' \ge 0 $, $s \in { S}_\phi$ with $|s t^{-1/{\hat n}}| > R$
\begin{equation}
\label{eq:eqintbound}
\left | \int_0^1 s^{l'} \lambda e^{-t 
\mathcal{P} \left ( -s t^{-1/{\hat n}} \right ) \lambda (1 - \tau) } d\tau \right | 
\le \lambda \left \{ \frac{1-e^{-t \lambda \Re \mathcal{P} (-s t^{-1/{\hat n}} )}}{t 
\lambda \Re P \left (-s t^{-1/{\hat n}} \right ) } \right \} s^{l'} \le C(\phi, T) 
t^{l'/{\hat n} - l'/n} 
\end{equation}
The definition of $C(\phi, T)$ implies that for $l' \ge 0$, $s \in {
  S}_\phi$ with $|s t^{-1/{\hat n}}| \le R$ we have
\begin{equation}
\label{eq:eqintbound2}
\left | \int_0^1 s^{l'} \lambda e^{-t 
\mathcal{P} \left ( -s t^{-1/{\hat n}} \right ) \lambda (1 - \tau) } d\tau \right | 
\le C(\phi, T) t^{l'/{\hat n} - l'/n} 
\end{equation}
Setting $l' = \sum j q_{l, j}$, using (\ref{eq:eqintbound}) and (\ref{eq:eqintbound2}),
we find
after time integration
\begin{multline}\label{nn1}
\Bigg\| \int_0^1 \lambda e^{-t \mathcal{P} (-s t^{-1/{\hat n}}) \lambda (1-\tau)}
{ {\bf B}}_{{\bf q}, {\bf k} } * { {\bf F}}^{*{\bf k}}
* \sideset{^*}{}\prod_{l=1}^{m}
\sideset{^*}{}\prod_{j=1}^{n} \left (
s^j { F}_l \right )^{*q_{l,j}}
(s, \lambda\tau; t) d\tau
\|_\nu  \\
\le t^{l'/{\hat n} - l'/n}
C(\phi, T) c_{{\bf q}, {\bf k}} (\nu, t) \| { {\bf F}} \Big\|_\nu^{
| {\bf q} |  + | {\bf k} | }
\end{multline}
Using (\ref{eqcond}), (\ref{eq:eqhatF}), (\ref{eq:eqfirst}) and (\ref{nn1}) , it follows
that ${ {\cal N}}$ maps a ball of radius $2 \|{\bf { F}_0} \|_0 $ into
itself.  Using Lemma \ref{L8}, (\ref{eq:eqintbound}) and (\ref{eq:eqintbound2}),
we obtain
\begin{multline*}
\Bigg\| \int_0^1
\lambda { {\bf B}}_{{\bf q}, {\bf k} } *
\Bigg \{ ({ {\bf F}}+{ {\bf h}})^{*{\bf k}}
* \sideset{^*}{}\prod_{l=1}^{m}
\sideset{^*}{}\prod_{j=1}^{n} \left (
s^j [{ F}_l + { h}_l ]  \right )^{*q_{l,j}}
\\-
{ {\bf F}}^{*{\bf k}}
*\sideset{^*}{}\prod_{l=1}^{m}
\sideset{^*}{}\prod_{j=1}^{n} \left (
s^j { F}_l \right )^{*q_{l,j}} \Bigg \}
(s, \lambda\tau; t) e^{-t \mathcal{P} (-s t^{1/{\hat n}} ) \lambda (1-\tau)} d\tau
\Bigg\|_\nu
\\
\le t^{l'/{\hat n} - l'/n}
C(\phi, T) (| {\bf q} |  + | {\bf k} |  )
c_{{\bf q}, {\bf k}} (\nu, t)
\left ( \| { {\bf h}} \|_\nu + \| { {\bf F}} \|_\nu \right )^{
| {\bf q} | + | {\bf k} | - 1 }
\| { {\bf h}} \|_\nu
\end{multline*}
where $l' = \sum j q_{l, j}$ from which the conclusion using
(\ref{defck}) and (\ref{eqcond}).
\end{proof}

\bigskip

\noindent{\bf Behavior of  $^s\!\mathbf { F}$ near $s=0$}

In the following proposition, we denote by 
$^s\!\mathbf F$ the solution $\bf F$ of Lemma 
\ref{eq:eqlemmasmall}. 

\begin{Proposition}
  \label{AsymptAts=0} For small $s$ we have  
$$|^s\!\mathbf { F}|\le \sum_{j=1}^{J_r} |s|^{-\omega_j-1} t^{(1+\omega_j)/{\hat n}}
\mathfrak{h}^+_{n'_j} (|s|^{\beta_1}, t^{\gamma_2-\beta_2/{\hat n}} |s|^{\beta_2}, ...
t^{\gamma_K-\beta_K/{\hat n}} |s|^{\beta_K} ) 
$$
\end{Proposition}
\begin{proof}
   The proof is similar to that of Proposition~\ref{AsymptAtp=0}, using (\ref{eq:eqF0b}), (\ref{eq:eqr}) and
(\ref{eq:eqn'}). 
  $^s\!{\bf { F}}$ to (\ref{eq:eqhatF}) solves a
  linear equation 
\begin{equation}
\label{29.2.0}
^s\!{\bf { F}}=\mathcal{G}\, (^s\!{\bf { F}})+{\bf { F}}_0\ \ \mbox{\rm or }\ \
^s\!{\bf { F}}=(1-\mathcal{G})^{-1}{\bf { F}}_0
\end{equation}
with ${\cal G}$ very
similar to that given in \S4.
\end{proof}
\smallskip

\noindent{\bf End of proof of Theorem \ref{Tasympt}} (i) The proof is
a direct application of Lemma \ref{eq:eqlemmasmall} and Proposition
\ref{AsymptAts=0}.  Using (\ref{eqlaplace}) and properties of Laplace
transform, (\ref{fsbound}) follows for large $|\zeta|$, in the sector
$\arg \zeta \in \left ( -\frac{\pi}{2} - \phi, \frac{\pi}{2} + \phi
\right )$.

\subsection{Proof of Theorem \ref{Tasympt} (ii)}
An important difference is that infinite sums appear in some estimates. 
Analyticity  of the functions  $\mathfrak{a}$ and the estimate
$$
\| \mathcal{L}^{-1} y^{-\alpha} \|_\nu = \left\| \frac{p^{\alpha
      -1}}{\Gamma (\alpha) } \right\|_\nu \le C (1+\alpha^2) \nu^{-\alpha +1} , $$
for $\nu > 1$ with $C$ is independent of $\alpha$  and $\nu$, 
show convergence of the
corresponding series.  Also, the proof of Lemma \ref{eq:eqlemmasmall}
holds if the following norm was used instead:
$$
\| { F} \|_\nu^u = \sup_{0 \le \lambda \le 1, |t| \le T, s \in
  \mathcal{ S}_\phi} (1+|s|^2) e^{-\nu |s|} |{ F} (s,\lambda;  t) | $$
since for ${\hat n} =n $, $\Re t \mathcal{P} (-st^{-1/n}) = \Re s^n $,
is independent of $t$ in the exponent in (\ref{eq:eqhatF}).  
To show analyticity, we let $\hat{G}(s,\lambda;\theta)=\theta^{-(1+\omega_1)/(n\omega)}\hat{F}(s,\lambda;\theta^{1/\omega})$; then 
$\hat{G}$ satisfies an equation of the form
$$\hat{G}=\mathcal{N}_1(\hat{G})$$
where the conditions in Setting 2
and the choice of $\omega$ are such that $\mathcal{N}_1$, as it is
seen after straightforward algebra, manifestly preserves analyticity in
$\theta$. Using (\ref{eqlaplace}), analyticity of $t^{-\omega_1/n} { f}
(\zeta, t)$ in $t^{\omega}$ follows provided $|\zeta|$ is large enough
(depending on $T$).
\subsection{Proof of Corollary~\ref{mhdym}}
Substitution gives for $f(x, t)$, defined by (\ref{eqHf}), an equation
of the form (\ref{1.a}), with $m=1$, $d=1$. Then in (\ref{3}), $\bf k$
is scalar. The vector $\bf q$ is $3$ dimensional, indexed by
$(l,j),\,l=1, j=1,2,3$.  The nonlinearity is
  quartic and the equation is linear in the derivatives of $f$, thus
  the only nonzero values of $b_{\mb q,k}$ are when $\bf q$ is $\bf 0$
  (and $k=1,...,4$) or a unit vector $\hat{\bf e}_i\in\RR^3$ (and
  $k=0,...,3$).  Further, it is found that
$$J_r=1, K=2, \omega_1 = \frac{5}{3} = \beta , \gamma_1=\gamma_2=1,
\beta_1 = 3, \beta_2 = \frac{2}{3}, {\hat n} = 3$$
and in
(\ref{eq:eqA2}) we have
\begin{equation}
\label{valpha}
\alpha_{{\bf 0}, 1} = \frac{4}{3} , ~ 
\alpha_{{\bf 0}, 2} = -1 ,~
\alpha_{\hat{\bf e}_1, 1} = 2 ,~ 
\alpha_{\hat{\bf e}_2, 1} = 1 , ~
\alpha_{\hat{\bf e}_3, 1} =  0
\end{equation}
This is sufficient to check that Theorem \ref{Tasympt} applies.

Since $|z| t^{-2/9}$ large 
corresponds to
$ |\zeta| = |x| t^{-1/3}$ large, and 
$\arg z \in \left (- \frac{4}{9} \pi, \frac{4}{9} \pi \right ) $
corresponds to $\arg \zeta \in \left ( -\frac{2}{3} \pi , \frac{2}{3} \pi
\right ) $, Theorem \ref{Tasympt}
implies that for any $\phi \in (0, \frac{\pi}{6} )$ for large 
$x \in \mathcal{D}_\phi $ and large $\zeta = x/t^{1/3}$ we have
$$
|f(x, t) | = O\left (|x|^{5/3} \mathfrak{h}_{(N+1)} (t |x|^{-3}, t |x|^{-2/3}
\right ) = O\left (|x|^{5/3} t^{N+1} \mathfrak{h}_{(N+1)} (|x|^{-3}, |x|^{-2/3}
\right ) $$
Changing variables, this implies
\begin{multline*}
x (z)^{-2} f (x(z,t), t) = 
O \left ( t^{N+1} |z|^{-\frac{1}{2}}  \mathfrak{h}_{(N+1)} (|z|^{-\frac{9}{2}}, |z|^{-1} \right ) 
\\= o \left ( t^{N} |z|^{-\frac{1}{2}}  \mathfrak{h}_{(N)} (|z|^{-\frac{9}{2}}, |z|^{-1} \right )
\end{multline*}
as needed for asymptoticity.
The convergence in the series representation in $t^{7/9}$ follows from
Theorem \ref{Tasympt} (ii).  It is seen from (\ref{eqrational}) that
all the exponents of $t$ are integer multiples of $\frac{7}{9}$.$\Box$

\begin{Note}  Large $\zeta$ includes part of the region where
  Theorems \ref{T1} and \ref{TrB} imply Borel summability of the
  expansion in inverse powers of $z$. Together, the results provide
  uniform control of the solution.
\end{Note}

\section{Appendix}
\subsection{Asymptotic behavior: further comments}\label{Asympts}

In the assumptions of Theorem~\ref{TrB}, by the remark following it,
formal series solutions to the initial value problem are asymptotic to
the actual unique solution.  The discussion below addresses the issue
of deriving this series, or, when less regularity is provided and only
the first few terms of the expansion exist, how to show their
asymptoticity.

\z {\em Heuristic calculation}. Assuming algebraic behavior of $\bf f$
in our assumptions on the nonlinearity, it is seen that the most
important terms for large ${\bf x}$ (giving the ``dominant balance'')
are ${\bf f}_t $, $\mathcal{P}_0 \bf f$, coming from the constant part
of $\mathcal{P}$, and ${\bf r} ({\bf x},t)$.  This suggests that, to
leading order, ${\bf f} ({\bf x}, t) \sim {\bf f}_I ({\bf x})
+\int_0^t e^{-\mathcal{P}_0 (t-\tau)} {\bf r} ({\bf x}, \tau) d\tau $.
If we substitute 
\begin{equation}
  \label{subs1}
  {\bf f} ({\bf x}, t) = \mathbf{A}_1 (t){\bf x}^{-\alpha_r
  \bf 1} + {\bf {\tilde f}}
\end{equation}
into (\ref{1.a}), ${\bf {\tilde f}}$ will generally satisfy an
equation of the form (\ref{1.a}), for an {\em increased} value of
$\alpha_{r}$; if the process can be iterated, as is the case in the
examples in \cite{CPAM}, it generates a formal series solution.

To obtain rigorous estimates, one writes the equation for ${\bf
  {\tilde f}}$ defined in (\ref{subs1}) and applies Theorem~\ref{T1}
to show ${\bf {\tilde f}}=o( {\bf x}^{-\alpha_{r}\bf 1} )$.  If the
coefficients of the equation allow it, this procedure can be repeated
to obtain more asymptotic terms for ${\mathbf{f}}$. This is the case
for instance in the assumptions of Theorem~\ref{TrB}, where a complete
series is obtained, which is furthermore Borel summable to $\bf f$.

The discussion also shows that the assumption $\alpha_r \ge 1$ can be
often be circumvented by subtracting the higher powers of $\bf x$ from
$\bf f$.

\subsection{Simple examples of Borel regularization}\label{illustr}
In this section we discuss informally and using rather trivial
examples, the regularizing features of Borel summation. An excellent
account of \'Ecalle's modern theory of generalized summability is
found in \cite{EcalleNato};  see \cite{TOP} as well.  
Many interesting results, using more
classical tools can be found in \cite{Balser}.

Singular perturbations give rise to nonanalytic behavior and divergent
series. Infinity is an irregular singular point of the ODE $f'-f=1/x$,
and the formal power series solution $\tilde{f}=\sum_{k=0}^\infty
(-1)^k k!x^{-k-1}$ diverges. In the context of PDEs, the solution $h$
of the heat equation $h_{t}-h_{xx}=0$ with $h(0,x)$ real-analytic but
not entire, has a factorially divergent expansion in {\em small $t$},
the recurrence relation for the terms of which is $kH_{k}=H_{k-1}''$.

The {\em Borel transform} of a series, is by definition its term-wise
inverse Laplace transform, which improves convergence since
$\mathcal{L}^{-1} x^{-k-1}=p^k/k!$.  If the Borel transformed of a
series converges to a function which can be continued analytically
along $\RR^+$ {\em and} is exponentially bounded, then its Laplace
transform is by definition 
the {\em Borel sum} of the series. Since on a formal
level Borel summation is $\mathcal{L}\mathcal{L}^{-1}$, the identity,
it can be shown to be an extended isomorphism between series and
functions; in particular, the Borel sum of $\tilde{f}$ above,
$\mathcal{L}(1+p)^{-1}$ is an actual solution of the equation. Another
way to view this situation is that Borel transform maps singular
problems into more regular ones. The Borel transform of the ODE
discussed is $(p+1)\mathcal L^{-1} f+1=0$. The inverse Laplace
transform of $h_t=h_{xx}$ in $1/t$ is
$\hat{h}_{xx}-p\hat{h}_{pp}-\frac{3}{2}\hat{h}_p=0$ which becomes
regular, $u_{xx}-u_{zz}=0$ by taking $\hat{h}(p,x)$
$=p^{-1/2}u(2 p^{1/2}, x)$, $z=2 p^{1/2}$. 

It is in its latter role, of a regularizing tool, that we
use Borel summation in PDEs.

\subsection{Derivation of equation (\ref{1.a}) from (\ref{1})}\label{D11}
 We define an $m$-dimensional vector
${\bf f}$ by ordering the set $\left \{ \partial_{\bf x}^{\bf j} {\bf
    u}: 0 \le |{\bf j} | <n \right\}$.  It is convenient to introduce
${\bf{\hat g}}_2 ({\bf x}, t, {\bf f})$ so that
$$ \sum_{|{\bf J}| = n} {\bf g_{2,J}} \left ( {\bf x}, t, 
\{ \partial_{\bf x}^{\bf j} {\bf u}  \}_{|{\bf j}|\le n-1} 
\right ) \partial_{\bf x}^{\bf J}  {\bf u} 
= -\sum_i 
{\bf {\hat g}_{2,i}} ({\bf x}, t, {\bf f}) \partial_{x_i} {\bf f} 
$$
So, for showing that (\ref{1}) implies
(\ref{1.a}) it is enough 
to show that for $1 \le n' \le n$, for $|{\bf J'} |= n'-1$,
$$\partial_{\bf x}^{\bf J'} 
\left [ {\bf g}_1 ({\bf x}, t, {\bf f}) + \sum_i 
{\bf {\hat g}}_{2_i} ({\bf x}, t, {\bf f} )  
\partial_{x_i} {\bf f} \right ] $$ 
is of the form on the right hand side of (\ref{1.a}). We do so in three
steps.

\begin{Lemma}
\label{A0}
Consider for  $k \ge 1 $,   
\begin{equation}
\label{A0.1}
{\bf E}({\bf x}, t) = 
{\sum_{{\bf q} \succeq 0 }}^\ddagger {\bf b_{\bf
q}} ({\bf x},t,{\bf f}) 
\prod_{\{m;k\}}
\left (\partial_{\bf x}^{\bf j} f_l
\right )^{q_{l,{\bf j}}} 
\end{equation}
where $\{m;k\}$ denotes the set
$\{(l,\mb j):1\le l\le m;1\le |\mb j |\le k\}$, and 
$\ddagger$ means summation over ${\bf q}$ with the restriction
\begin{equation}
\label{A0.2}
\sum_{\{m;k\}}
| {\bf j} | q_{l, {\bf j}} \le k 
\end{equation}   
Then, for $i =1, 2..,d $, 
$\partial_{x_i} {\bf E} ({\bf x}, t) $ has the same form
as (\ref{A0.1}) with restriction (\ref{A0.2}), provided
$k$ is replaced by $k+1$.    
\end{Lemma}
\begin{proof}
The proof is straightforward, keeping track of the number of derivatives and the powers involved: 
note that 
$$\partial_{x_i} {\bf E} ({\bf x}, t, {\bf f})    
= \sum_{{\bf q} \succeq 0 } \left ( \sum_{l=1}^m 
{\frac {\partial}{\partial f_l}} {\bf b_{\bf
q}} ({\bf x},t,{\bf f}) \partial_{x_i} f_l + \partial_{x_i} 
{\bf b}_{{\bf q}} ({\bf x}, t, {\bf f}) \right ) 
\prod_{\{m;k\}}
\left (\partial_{\bf x}^{\bf j} f_l
\right )^{q_{l,{\bf j}}} 
$$
$$ +  
\sum_{{\bf q} \succeq 0 } 
{\bf b}_{{\bf q}} ({\bf x}, t, {\bf f}) ) 
\sum_{l'=1}^m 
\sum_{| {\bf j'}|=1}^{k} 
q_{l',{\bf j'}} \left ( \partial_{\bf x}^{\bf j'} 
f_{l'} \right )^{q_{l',{\bf j'}}-1} 
\partial_{x_i} (\partial_{\bf x}^{\bf j'} f_{l'} )  
\mathop{{\prod}^\dagger}_{\{m;k\}} 
\left (\partial_x^j f_l
\right )^{q_{l,{\bf j}}} 
$$
where $\prod^{\dagger}$ indicates that the term $l=l',
{\bf j}={\bf j'}$ is missing from the product.  Manifestly, this
  is of the form (\ref{A0.1}) with a suitable redefinition of ${\bf
  b}_{\bf q}$ and with the product of the number of derivatives times
the power totaling at most
$$ | {\bf j'}| + 1 + 
|{\bf j'}| (q_{l', {\bf j}'} - 1)  
+ {\sum_{\{m; k\}}}^\dagger 
| {\bf j} | q_{l, {\bf j}} 
= 1   
+ \sum_{\{m; k\}} 
|{\bf j}| q_{l, {\bf j}} \le k+1 $$
Hence restriction (\ref{A0.2}) holds, now with $k+1$ instead of $k$.
  \end{proof}

\begin{Lemma} 
\label{A1}
For any $n'\ge 1$, and any ${\bf J'}$ with 
$| {\bf J'} | = n'-1$,
\begin{equation}
\label{A.1}
\partial_{\bf x}^{\bf J'} {\bf g_1} (y, t, {\bf f} (y, t))    
= {\sum_{{\bf q} \succeq 0 }}^\ddagger {\bf b_{\bf
q}} ({\bf x},t,{\bf f}) 
\prod_{\{m; n'-1\}}
\left (\partial_{\bf x}^{\bf j} f_l
\right )^{q_{l,{\bf j}}}
\end{equation}

\z for some ${\bf b}_{\bf q}$, depending on $n^\prime$, $\bf g_1$,
and its first $n'-1$ derivatives, and where $\sum^\ddagger$ means the sum over ${\bf q} $
with the further restriction
$$
\sum_{\{m;n'-1\}}
| {\bf j} | q_{l, {\bf j}} \le n'-1 
$$
\end{Lemma}
\begin{proof}
The proof is by induction. We have, with obvious notation,
$$\partial_{x_i} {\bf g}_1 ({\bf x}, t, {\bf f} ({\bf x}, t)) = 
{\bf g}_{1,x_i} + 
{\bf g}_{1,\mathbf f}\cdot \partial_{x_i} \mathbf f $$

\z which  is of the form
(\ref{A.1}).  Assume (\ref{A.1}) holds for $n' = k \ge 1$, i.e.
for all ${\bf J}'$ satisfying $|{\bf J'} | = k-1$,
$$
\partial_{\bf x}^{\bf J'} {\bf g_1} ({\bf x}, t, {\bf f}) 
= {\sum_{{\bf q} \succeq 0 }}^\ddagger 
{\bf b_{\bf q}}
({\bf x} ,t,{\bf f}) 
\prod_{\{m; k-1\}} 
\left (\partial_{\bf x}^{\bf j} f_l
\right )^{q_{l,{\bf j}}}
$$
Taking a $x_i$ derivative, and applying Lemma \ref{A0}, $
\partial_{\bf x}^{\bf J} {\bf g_1} (y, t, {\bf f}) $ for $|{\bf J} | =
k $ will have the form above, with $k-1$ replaced by $k$ and with
restriction
$$ \sum_{\{m;k\}}
| {\bf j} | q_{l, {\bf j}} \le k 
$$
Thus, (\ref{A.1}) holds for $n'=k+1$, with a different ${\bf b}$. 
The induction step is proved.
\end{proof}

\begin{Lemma} 
\label{A2}
For $n'=1,2,...,n$, and any ${\bf J}$ with $|{\bf J}| = n'-1$ we
have
\begin{equation}
\label{A.11}
\partial_{\bf x}^{\bf J} \left [ {\bf {\hat g}_{2,i'}} 
({\bf x}, t, {\bf f}) \partial_{\bf x_{i'}} {\bf f}
\right ] = 
{\sum_{{\bf q} \succeq 0 }}^\ddagger 
{\bf b_{\bf q}} ({\bf x},t,{\bf
f}) 
\prod_{\{m; n'\}}
\left (\partial_{\bf x}^{\bf j} f_l \right
)^{q_{l,{\bf j}}}\end{equation} 
for some ${\bf b}_{{\bf q}}$,
depending on $n^\prime$, $\mathbf g_2$ 
and its first $n'-1$ derivatives, where $\sum_{{\bf q}\succeq 0}^{\ddagger}$
denotes summation with the restriction
\begin{equation}
\label{A.11.1}
\sum_{\{m;n'\}}
|{\bf j} | q_{l, {\bf j}} \le n^\prime
\end{equation} 
\end{Lemma}

\begin{proof}
  Clearly (\ref{A.11}) with restriction (\ref{A.11.1}) holds for
  $n'=1$. Suppose it holds for $n'=k$. Then we
  note that if $| {\bf J} | = k+1$, then there exists some
  index $1 \le i \le d$ and some ${\bf J'}$, with $| {\bf J'} | = k$
  so that $\partial_{\bf x}^{\bf J} = \partial_{x_i} [ \partial_{\bf
    x}^{\bf J'} ]$; hence applying Lemma \ref{A0}, we obtain 
   (\ref{A.11}) and (\ref{A.11.1}) for $n'=(k+1)$.
\end{proof}

\subsection{Some useful inequalities.}\label{Proof37}
\begin{enumerate}
\item We start with a simple inequality for $\alpha>1$ and $\mu >0$:
\begin{equation}
  \label{ine1}
  (1+\mu ^\alpha)\int_0^1s^{\alpha-1}e^{-\mu s}ds\le 2\Gamma(\alpha)
\end{equation}
This is clear for $\mu \le1$, while for $\mu >1$ we write $(1+\mu ^\alpha)\le 2
\mu ^\alpha$ and note that $\int_0^\infty
s^{\alpha-1}e^{-\mu s}ds=\mu ^{-\alpha}\Gamma(\alpha)$.  
\item  For $\alpha > 0 $, $\mu  > 0$, $\sigma=0,1$, $\nu >  2$ and 
$m\in\NN$, 
\begin{equation}
  \label{minilemma}
  \mu^\alpha \nu^\alpha\int_0^1 \frac{e^{-\nu \mu  [1 -(1-s)^m]}}{[1 + \mu ^2 (1-s)^2]^\sigma} s^{\alpha - 1} ds 
\le 8(2^{\alpha}+1) \Gamma (\alpha)  [1+\mu ^2 ]^{-\sigma} 
\end{equation}
where $C(m)$ is independent of $\mu $, $\alpha$ and $\nu$.
Indeed,  the integral is bounded by
\begin{multline*}
\Big(\int_0^{\frac{1}{2}}du+\int_{\frac{1}{2}}^1du\Big)\frac{e^{-\mu\nu s}s^{\alpha-1} ds}{[1 + \mu ^2 (1-s)^2]^\sigma}\le \frac{1}{(1+\mu^2/4)^\sigma}\int_0^1 e^{-\mu\nu s}s^{\alpha-1}ds
\\
+\max_{s\in[1/2,1]}\frac{e^{-\mu\nu s}}{[1 + \mu ^2 (1-s)^2]^\sigma}\int_0^1s^{\alpha-1}ds\le \frac{2\Gamma(\alpha)(\mu\nu)^{-\alpha}}{(1+\mu^2/4)^\sigma}+\frac{e^{-\mu\nu/2}}{\alpha(1+\mu^2/4)^\sigma}\\
\le \frac{2\Gamma(\alpha)(\mu\nu)^{-\alpha}}{(1+\mu^2/4)^\sigma}+ \frac{2^{\alpha+1}\Gamma(\alpha)(\mu\nu)^{-\alpha}}{(1+\mu^2/4)^\sigma}\sup_{\alpha\in\RR^+}
\sup_{\mu\nu\in\RR^+}\frac{(\mu\nu)^\alpha e^{-\mu\nu/2}}{2^{\alpha+1}\alpha\Gamma(\alpha)}\\\le
\frac{2\Gamma(\alpha)(\mu\nu)^{-\alpha}}{(1+\mu^2/4)^\sigma}+ \frac{2^{\alpha+1}\Gamma(\alpha)(\mu\nu)^{-\alpha}}{(1+\mu^2/4)^\sigma}
\end{multline*}

\item For $n>1$ the function
\begin{equation*}
  \label{bdR}
  (1+\mu )e^{-\mu }\int_0^1e^{\mu [u^n+(1-u)^n]}du
\end{equation*}
is bounded in $\RR^+$, as it can be checked applying Watson's lemma
for large $\mu $ and noting its continuity on $[0,\infty)$.  Thus, for
some constant $C$ and $\nu>1$ we have
\begin{equation}
  \label{pd}
  \int_0^{|p|}e^{\nu |s|^n+\nu |p-s|^n}ds\le \frac{C|p|}{1+|p|^n}e^{\nu |p|^n}
\end{equation}

\item  We have $|\mathbf{p}^{\mathbf{k}}|\le \max_{i\le d}
  |p|_i^{|\mathbf{k}|}\le \sum_{i\le
    d}|p_i|^{|\mathbf{k}|}$ and thus for some constant $C$ and all $j\le m$ we have
\begin{equation}
  \label{opA}
  |\mathcal{P}_j(-\mathbf{p})|\le
C\sum_i(1+|p_i|^n)
\end{equation}
Also, for some $C_2>0$, $| \mathcal{P}_j (-{\bf p})|\le C_2\sum_{i}(1+|
p_i|+| p_i^n|)=:C_2 (d+q)$ and thus, for $\nu>C_2+1$ we have, for $0\le l'\le n$,
\begin{multline}
    \label{int2}
   |{\bf p}|^{l'} \int_{0}^t e^{| \mathcal{P}_j (-{\bf p})|
      (t-\tau)}e^{\nu(\tau+1)q}d\tau\le  
|{\bf p}|^{l'}e^{q\nu+C_2 t d }\int_0^te^{(\nu-C_2)q\tau}  d\tau\\
\le T^{1-l'/n} e^{\nu q (t+1) + C_2 t d} 
 \frac{|{\bf p}|^{l'}}{[(\nu - C_2) q]^{l'/n}} \sup_{\gamma > 0} 
\frac{1 - e^{-\gamma}}{\gamma^{1-l'/n}}
\le \frac{C_3 (T)}{(\nu - C_2)^{l'/n}} 
e^{\nu q (t+1) + C_2 t d} 
\end{multline}
\end{enumerate}

\subsection{Modified estimates for Lemma~\ref{combined}.}\label{PL31}

>From (\ref{pd}) it follows that for a constant $C$ independent of $
\boldsymbol{\Psi},\boldsymbol{\Phi}$ we have
\begin{equation}\label{mod}
|\boldsymbol{\Psi}*\boldsymbol{\Phi}|\le Ce^{\nu (t+1)\sum_i(|p_l|+|p_l|^n)} \|\boldsymbol{\Psi}\|_{\nu n} \|\boldsymbol{\Phi}\|_{\nu n}
\end{equation}
In particular $\mathfrak{B}(\nu,n,\mathcal{S})$ is a Banach algebra.
For the equivalent of Lemma~\ref{L4}, we use the following bounds.
\begin{multline}\label{new3}
I=\int_0^{|p_1|} s^{\alpha - 1} e^{- \nu (t+1) [|p_1|^n - (|p_1|-s)^n]}  
e^{-\nu (t+1) s} ds\le \int_0^{|p_1|} s^{\alpha - 1} e^{-\nu (t+1) s} ds 
 \\\le \frac{\nu^{-\alpha}}{\Gamma(\alpha)(t+1)^{\alpha}}\\
{\text{and }}\  I
\le |p_1|^\alpha \int_0^{1} s^{\alpha - 1} e^{-\nu (t+1) |p_1|^n 
[ 1- (1-s)^n ] } ds \le 
C \frac{2^\alpha \Gamma (\alpha) |p_1|^\alpha}{[\nu (t+1) |p_1|^n]^\alpha}    
\end{multline}
where we used (\ref{minilemma}) for $\sigma = 0$.  From (\ref{new3})
it is clear that
\begin{equation}\label{181} 
\|{\bf H} * F_j \|_{\nu n} \le \big \| |{\bf H}| * |F_j| \big \|_{\nu n}
\le C[\Gamma (\alpha)]^{d}c^\alpha (\nu(t+1))^{-d \alpha} \| {\bf F} \|_{\nu n}
\end{equation} In Lemma \ref{L4.7.9}, we get instead
$$
\Big | |{\bf F}| * |{\bf G}| \Big | \le e^{\nu(t+1) \sum_i
  (|p_i|+|p_i|^n) } \| {\bf F} \|_{\nu n} \| {\bf G} \|_{\nu n} $$
Very similar changes are made in in Lemma~\ref{L4.9.5},
Corollary~\ref{C5}, and in Lemma~\ref{L6} where in the proof we use
(\ref{int2}) instead of (\ref{eq:28.5}).  Definition~\ref{D4},
Lemma~\ref{L7} and Definition~\ref{D41} do not change.
Lemma~\ref{L7.5}, Lemma~\ref{L8} change in the same way as above. In
Lemma~\ref{L8.9} we use again (\ref{int2}) instead of (\ref{eq:28.5})
to make corresponding changes. Finally, in Lemma~\ref{L10},
$\nu/4$ changes to $\nu/4/c$.

\section{Acknowledgments} The authors are very grateful to B L J Braaksma for
  a careful reading of the manuscript and many useful suggestions. The
authors are indebted to R D Costin for valuable suggestions. One of the
authors also benefited from discussions with B. Sandstede. Work supported by
NSF Grants DMS-0100495, DMS-0074924, DMS-0103829.  Travel support by the Math
Research Institute of the Ohio State University is also gratefully
acknowledged.

\vfill \eject

\newpage
\figure 
\ifx\pdftexversion\undefined
 
\else
  $ $ \vskip -8cm
\fi
\includegraphics{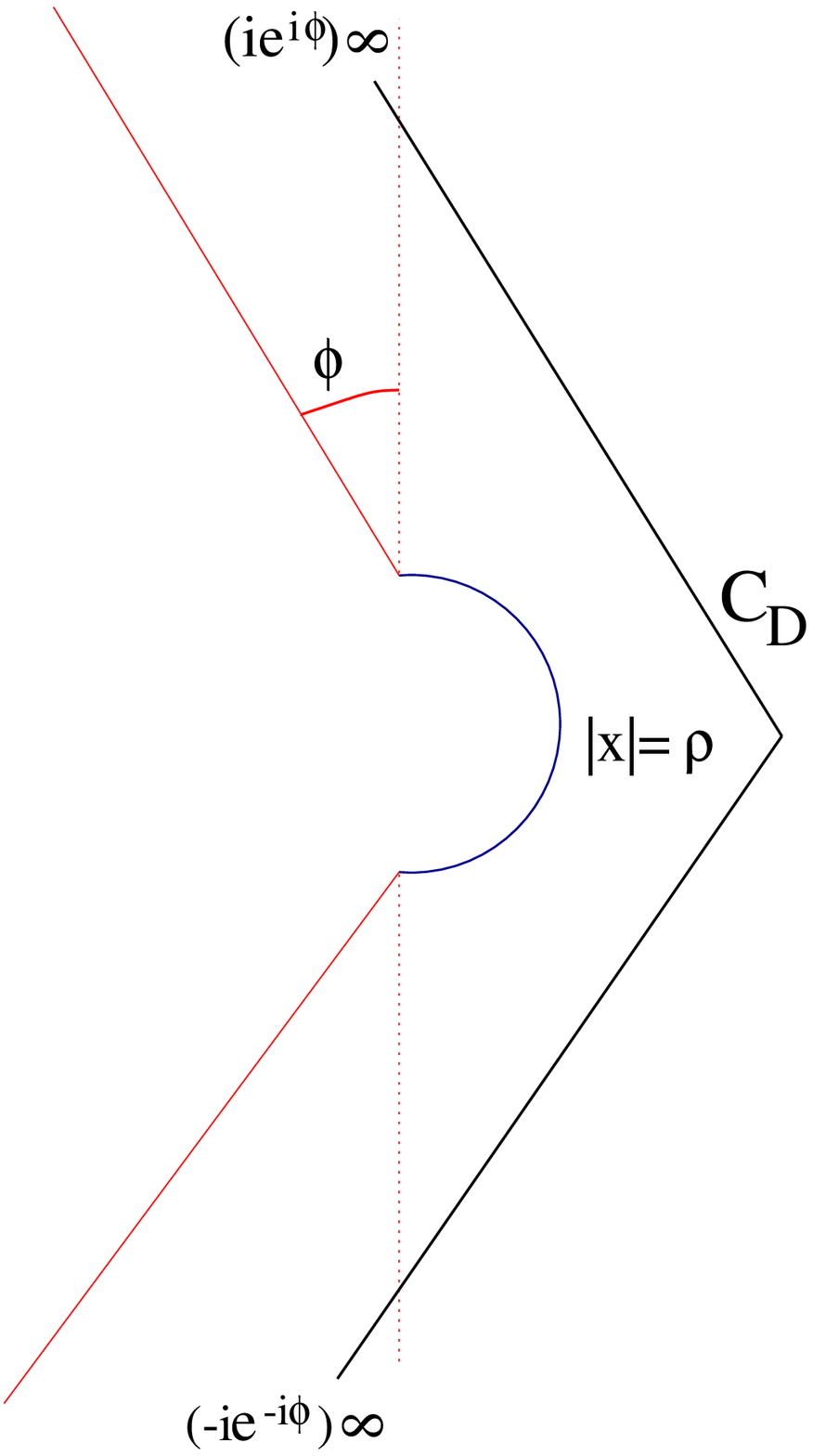}
\caption{Contour $C_D$ in the $(\mb p)_i-$plane.}
\endfigure

\end{document}